\documentclass{article}

\usepackage{fancyhdr}
\usepackage[T1]{fontenc}
\usepackage[english]{babel}
\usepackage[all]{xy}
\usepackage{amsfonts}
\usepackage{amsmath}
\usepackage{amssymb}
\usepackage{latexsym}
\usepackage{mathrsfs} 
\usepackage{graphicx} 
\usepackage[toc,page]{appendix}

\newtheorem{thm}{\bf Theorem}[section]
\newtheorem{cor}[thm]{\bf Corollary}

\newtheorem{lem}[thm]{\bf Lemma}
\newtheorem{prop}[thm]{\bf Proposition}
\newtheorem{defn}[thm]{\bf Definition}

\newtheorem{rem}[thm]{\bf Remark}

\newtheorem{exmp}[thm]{\bf Example}

\newcommand{\field}[1]{\mathbb{#1}}
\newcommand{\C}{\field{C}}
\newcommand{\R}{\field{R}}

\newcommand{\N }{\field{N}}

\def\proof{{\parindent0pt {\bf Proof.\ }}}

\def\pd{{\rm pd}}
\def\fd{{\rm fd}}
\def\id{{\rm id}}

\def\Gpd{{\rm Gpd}}

\def\Gid{{\rm Gid}}

\def\pd{{\rm pd}}
\def\Gpd{{\rm Gpd}}
\def\P{\mathscr{P}_C}

\def\Gpc{{\rm G_C\!-\!pd}}
\def\1Gpc{\rm G_{C_1}\!-\!pd}
\def\2Gpc{\rm G_{C_2}\!-\!pd}

\def\G1ci{{\rm G_{C_1}\!-\!id}}
\def\Gci2{{\rm G_{C_2}\!-\!id}}

\def\Im{{\rm Im}}
\def\Coker{{\rm Coker}}
\def\Ker{{\rm Ker}}

\def\Ext{{\rm Ext}}
\def\Tor{{\rm Tor}}
\def\Hom{{\rm Hom}}
\def\End{{\rm End}}
\def\Prod{{\rm Prod}}
\def\Add{{\rm Add}}
\def\add{{\rm add}}
\def\Proj{{\rm Proj}}
\def\Inj{{\rm Inj}}

\def\Ggldim{{\rm Ggldim}}

\newcommand{\cqfd}
{\hspace{1cm}
\rule{2mm}{2mm}%
\medbreak%
\par%
}

\begin{document}
\baselineskip=15pt

\title{Relative Gorenstein dimensions over triangular matrix rings}

\author{Driss Bennis$^1,a$ \hskip 2cm  Rachid El Maaouy$^1,b$ \\ \\ J. R. Garc\'{\i}a Rozas$^2,c$ \hskip 1,5cm Luis Oyonarte$^2,d$}

\date{}

\maketitle

\begin{center}
\small{1: CeReMaR Research Center, Faculty of Sciences, B.P. 1014, Mohammed V University in Rabat, Rabat, Morocco.
	
	\noindent a:   driss.bennis@fsr.um5.ac.ma; driss$\_$bennis@hotmail.com

	\noindent b:  rachid$\_$elmaaouy@um5.ac.ma; elmaaouy.rachid@gmail.com
	
	2: Departamento de  Matem\'{a}ticas, Universidad de Almer\'{i}a, 04071 Almer\'{i}a, Spain.
	
	\noindent  c:  jrgrozas@ual.es
	
	\noindent d:  oyonarte@ual.es}
\end{center}
\begin{abstract} 
	Let $A$ and $B$ be rings, $U$ a $(B,A)$-bimodule and  $T=\begin{pmatrix} A&0\\U&B \end{pmatrix}$ the triangular matrix ring. In this paper, several notions in relative Gorenstein algebra over  a triangular matrix ring are investigated.
 We first study  how to construct w-tilting (tilting, semidualizing) over $T$ using  the corresponding ones  over $A$ and $B$.  We show that when $U$ is relative (weakly) compatible  we are able to  describe the structure of $G_C$-projective modules over $T$. As an application, we study when a morphism in $T$-Mod has  a special $G_CP(T)$-precover and when the class $G_CP(T)$ is a special precovering class. In addition, we study the relative global dimension of $T$. In some cases,   we   show  that it can be computed from the  relative global dimensions of  $A$ and $B$.  We end the paper with a counterexample to a result that   characterizes when a $T$-module has a finite projective dimension.
\end{abstract}



\section{Introduction} \label{Sec:1}

Semidualizing modules were independently studied (under different names) by Foxby \cite{F}, Golod \cite{Go}, and Vasconcelos \cite{V} over a commutative Noetherian ring. Golod used these modules to  study $G_C$-dimension for finitely generated modules. Motivated (in part) by Enochs and Jenda's extensions of the classical G-dimension given in \cite{EJ}, Holm and J$\phi$rgensen, extended in \cite{HJ} this notion to arbitrary modules. After that, several generalizations of semidualizing and $G_C$-dimension have been made by several authors (\cite{Wh},\cite{LHX},\cite{ATY}).
 
As the authors mentioned in \cite{BGO1}, to study the Gorenstein projective modules and dimension relative to a semidualizing $(R,S)$-bimodule $C$, the condition $\End_S(C)\cong R$, seems to be too restrictive and in some cases unnecessary. So the authors introduced weakly Wakamatsu tilting as a weakly notion of semidualizing which made the theory of relative Gorenstein homological algebra wider and less restrictive but still consistent.  Weakly Wakamatsu tilting modules were subject of many publications which showed how important these modules could become in developing the theory of relative (Gorenstein) homological algebra (\cite{BGO1},\cite{BGO2},\cite{BDGO})

Let $A$ and $B$ be rings and $U$ be a $(B,A)$-bimodule. The ring $T=\begin{pmatrix}
A&0\\U&B
\end{pmatrix}$ is known as the formal triangular matrix ring with usual matrix addition
and multiplication. Such rings play an important role in the representation theory of algebras. The modules over such rings can be described in a very concrete fashion. So, formal triangular matrix rings and modules over them have proved to be a rich source of examples and counterexamples. Some important Gorenstein notions over formal triangular matrix rings have been studied by many authors (see \cite{Z, EIT, ZLW}). For example, Zhang \cite{Z} introduced compatible bimodules and explicitly described the Gorenstein projective modules over triangular matrix Artin algebra. Enochs, Izurdiaga, and Torrecillas \cite{EIT} characterized when a left module over a triangular matrix ring is Gorenstein projective or Gorenstein injective under the "Gorenstein regular" condition. Under the same condition, Zhu, Liu, and Wang \cite{ZLW} investigated Gorenstein homological dimensions of modules over triangular matrix rings. Mao \cite{M2} studied Gorenstein flat modules over $T$ (without the  "Gorenstein regular" condition) and gave an estimate of the weak global Gorenstein dimension of $T$.
 
The main objective of the present paper is to study relative Gorenstein homological notions (w-tilting, relative Gorenstein projective modules, relative Gorenstein projective dimensions and relative global projective dimension) over triangular matrix rings.
 
This article is organized as follows:

In Section 2, we give some preliminary results.

In Section 3, we study how to construct w-tilting (tilting, semidualizing) over $T$ using w-tilting (tilting, semidualizing) over $A$ and $B$ under the condition that $U$ is relative (weakly) compatible. We introduce (weakly) $C$-compatible $(B,A)$-bimodules  for a $T$-module $C$ (Definition \ref{C-compatible}). Given two w-tilting modules $_AC_1$ and $_BC_2$, we prove in Proposition \ref{when p preserves w-tilting} that $C=\begin{pmatrix}
C_1\\(U\otimes_AC_1)\oplus C_2
\end{pmatrix}$ is a w-tilting  $T$-module  when $U$ is weakly $C$-compatible. 

In Section 4, we first describe relative Gorenstein projective modules over $T$. Let $C=\begin{pmatrix}
C_1\\(U\otimes_AC_1)\oplus C_2
\end{pmatrix}$ be a $T$-module. We prove in Theorem \ref{structure of $G_C$-projective}, that if $U$ is $C$-compatible then a $T$-module $M=\begin{pmatrix}
M_1\\M_2
\end{pmatrix}_{\varphi^M}$ is $G_C$-projective if and only if $M_1$ is $G_{C_1}$-projective $A$-module, $\Coker\varphi^M$ is $G_{C_2}$-projective $B$-module and $\varphi^M:U\otimes_A M_1\to M_2$ is injective. As an application, we prove that the converse of Proposition \ref{when p preserves w-tilting} and refine in relative setting (Proposition \ref{when T is relative CM-free}), a result  of when $T$ is left (strongly) CM-free due to  Enochs, Izurdiaga, and Torrecillas in \cite{EIT}.
Also when $C$ is w-tilting, we characterize when a $T$-morphism is a special precover (see Proposition \ref{when a module has precover}). Then in Theorem \ref{the class of G_C-proj is a special precovering}, we prove that class of $G_C$-projective $T$-modules is special precovering if and only if so are the classes of $G_{C_1}$-projective $A$-modules and $G_{C_2}$-projective $B$-modules, respectively.

Finally, in Section 5, we give an estimate of $G_C$-projective dimension of a left $T$-module and the left $G_C$-projective global dimension of $T$. First, it is proven that, given a $T$-module $M=\begin{pmatrix}
M_1\\M_2
\end{pmatrix}_{\varphi^M}$, if $C=\textbf{p}(C_1,C_2):=\begin{pmatrix}
C_1\\(U\otimes_AC_1)\oplus C_2
\end{pmatrix}$ is w-tilting, $U$ is $C$-compatible and $$SG_{C_2}-PD(B):=sup\{\2Gpc(U\otimes_A G)\;|\; G\in G_{C_1}P(A)\}<\infty,$$ then 
 $$max\{\1Gpc(M_1),(\2Gpc(M_2))-(SG_{C_2}-PD(B))\}$$
 $$\leq \Gpc(M)\leq$$
 $$ max\{(\1Gpc(M_1))+(SG_{C_2}-PD(B))+1,\2Gpc(M_2)\}.$$
 As an application, we prove that, if $C=\textbf{p}(C_1,C_2)$ is w-tilting and $U$ is $C$-compatible then
 $$max\{G_{C_1}-PD(A),G_{C_2}-PD(B)\}$$
$$\leq G_C-PD(T)\leq $$
$$ max\{G_{C_1}-PD(A)+SG_{C_2}-PD(B)+1,G_{C_2}-PD(B)\}.$$  
Some cases when this estimation becomes an exact formula are also given. 

The authors in \cite{AS} establish a relationship between the projective  dimension of
modules over $T$ and modules over $A$ and $B$. Given  an integer $n\geq 0$  and $M=\begin{pmatrix}
M_1\\M_2
\end{pmatrix}_{\varphi^M}$  a $T$-module, they proved that $\pd_T(M)\leq n $ if and only if $\pd_A(M_1)\leq n$, $\pd_B(\overline{M}_2)\leq n$ and the map related to the $n$-th syzygy of $M$ is injective. We end the paper by giving a counterexample to this  result (Example \ref{counter example}).

\section{Preliminaries}
Throughout this paper, all rings will be associative (not necessarily commutative)
 with identity, and all modules will be, unless otherwise specified, unitary left modules.
For a ring $R$, We use $\Proj(R)$ (resp., $\Inj(R)$) to denote the class of all projective (resp.,
injective) $R$-modules. The category of all left $R$-modules will be denoted by $R$-Mod. For an $R$-module $C$, we use $\Add_R(C)$ to denote the class of all $R$-modules
which are isomorphic to direct summands of direct sums of copies of $C$, and
$\Prod_R(C)$ will denote the class of all $R$-modules which are isomorphic to direct
summands of direct products of copies of $C$.

Given a class of modules $\mathcal{F}$ (which will always be considered closed under isomorphisms), an $\mathcal{F}$-precover of $M\in R$-Mod is a morphism $\varphi:F\to M$ ($F\in \mathcal{F}$) such that $\Hom_R(F',\varphi)$ is surjective for every $F'\in\mathcal{F}$. If, in addition, any solution of the equation $\Hom_R(F,\varphi)(g)=\varphi$ is an automorphism of $F$, then $\varphi$ is said to be an $\mathcal{F}$-cover. The $\mathcal{F}$-precover $\varphi$ is said to be special if it is surjective and $\Ext^1(F,\ker \varphi)=0$ for every $F\in \mathcal{F}$. The class $\mathcal{F}$ is said to be special (pre)covering if every module has a special $\mathcal{F}$-(pre)cover.

Given the class $\mathcal{F}$, the class of all modules $N$ such that $\Ext^i_R(F,N)=0 $ for every $ F\in \mathcal{F}$ will be denoted by $\mathcal{F}^{\perp_i}$ (similarly, $^{\perp_i}\mathcal{F}=\{ N;\ \Ext^i_R(N,F)=0\ \forall F\in \mathcal{F}\}$). The right and left orthogonal classes $\mathcal{F}^{\perp}$ and $^{\perp}\mathcal{F}$ are defined as follows:
$$\mathcal{F}^{\perp}=\cap_{i\geq 1} \mathcal{F}^{\perp_i}\text{ and } ^{\perp}\mathcal{F}=\cap_{i\geq 1}\; ^{\perp_i}\mathcal{F}$$

 It is immediate to see that if $C$ is any module then $\Add_R(C)^{\perp}=\{C\}^{\perp}$ and ${^{\perp}\Prod_R(C)}={^{\perp}\{C\}}$.

Given a class of $R$-modules $\mathcal{F}$, an exact sequence of $R$-modules
$$\cdots\to X^1\to X^0\to X_0\to X_1\to\cdots$$
is called $\Hom_R(-,\mathcal{F})$-exact (resp., $\Hom_R(\mathcal{F},-)$-exact) if the functor $\Hom_R(-,F)$ (resp., $\Hom_R(F,-)$) leaves the sequence exact whenever $F\in \mathcal{F}$. If $\mathcal{F}=\{F\}$, we simply say $\Hom_R(-,F)$-exact. Similarly, we can define  $\mathcal{F}\otimes_R-$exact sequences when $\mathcal{F}$ is a class of right $R$-modules.

We now recall some concepts needed throughout the paper.

\begin{defn}
 \begin{enumerate}
 	\item (\cite[Definition 2.1]{HW}) A semidualizing bimodule is an $(R,S)$-bimodule $C$ satisfying the following
 	properties: 
	\begin{enumerate}
		\item $_RC$ and $C_S$ both admit a degreewise finite projective resolution in the corresponding module categories ($R$-Mod and Mod-$S$).
		\item $\Ext_R^{\geq 1}(C,C)=\Ext_S^{\geq 1}(C,C)=0.$
		\item The natural homothety maps $R\stackrel{_R\gamma}{\rightarrow}\Hom_S(C,C)$ and $S \stackrel{\gamma_S}{\rightarrow} \Hom_R(C,C)$ both are ring isomorphisms.
	\end{enumerate}

\item (\cite[Section 3]{W}) A Wakamatsu tilting module, simply tilting, is an $R$-module $C$  satisfying the following
properties: 
\begin{enumerate}
	 \item $_RC$ admits a degreewise finite projective resolution.
     \item $\Ext^{\geq1}_R(C,C)=0$
     \item   There exists a $\Hom_R(-,C)$-exact exact sequence of $R$-modules
     $$ \mathbf{X}=\ 0  \rightarrow R \rightarrow C^0 \rightarrow C^1 \rightarrow \cdots $$
     where  $C^i\in \add_R(C)$ for every $i\in \N $.
\end{enumerate}
 \end{enumerate}
\end{defn}
It was proved in \cite[Corollary 3.2]{W}, that an $(R,S)$-bimodule $C$ is semidualizing if and only if $_RC$ is tilting with $S=\End_R(C)$. So the following notion, which will be crucial in this paper, generalizes both concepts.
\begin{defn}[\cite{BGO1}, Definition 2.1] An $R$-module $C$ is weakly Wakamatsu tilting ($w$-tilting for short) if it has the following two properties:
\begin{enumerate}
\item  $\Ext_R^{i\geq 1}(C,C^{(I)})=0$ for every set $I.$
\item There exists a $\Hom_R (-,\Add_R(C))$-exact exact sequence of $R$-modules
 $$ \mathbf{X}=\ 0  \rightarrow R \rightarrow A^0 \rightarrow A^1 \rightarrow \cdots $$
 where  $A^i\in \Add_R(C)$ for every $i\in \N $.\medskip
\end{enumerate}
If $C$ satisfies $1.$ but perhaps not $2.$ then $C$ will be said to be $\Sigma$-self-orthogonal.
\end{defn}

\begin{defn}[\cite{BGO1}, Definition 2.2]
	Given any $C\in R$-Mod, an $R$-module $M$ is said to be G$_C$-projective if there exists a $\Hom_R (-,\Add_R(C))$-exact exact sequence  of $R$-modules
	$$ \mathbf{X}=\ \cdots \rightarrow P_1 \rightarrow P_0 \rightarrow A^0 \rightarrow A^1 \rightarrow \cdots $$
	where the $P_i's$ are all projective, $A^i\in \Add_R(C)$ for every $i\in \N $, $M \cong \Im(P_0 \rightarrow A^0)$.\medskip
	
	We use $G_CP(R)$ to denote the class of all G$_C$-projective $R$-modules.

\end{defn} 

It is immediate from the definitions that w-tilting modules can be characterized as follows.
\begin{lem}\label{charac of w-tilitng R-mod}
	An $R$-module $C$ is w-tilting if and only if both $C$ and $R$ are  $G_C$-projective modules.
\end{lem}

\medskip

Now we recall some facts about triangular matrix rings. Let $A$ and $B$ be rings and $U$ a $(B,A)$-bimodule. We shall denote by $T=\begin{pmatrix}
      A & 0\\ 
      U & B  
    \end{pmatrix}$ the generalized triangular matrix ring.
By \cite[Theorem 1.5]{G}, the category $T$-Mod of left $T$-modules is equivalent to the category $_T\Omega$ whose objects are triples $M=\begin{pmatrix}
      M_1 \\ 
      M_2   
    \end{pmatrix}_{\varphi^M}$, where $M_1\in A$-Mod, $M_2\in B$-Mod and $\varphi^M : U \otimes_A M_1\rightarrow M_2$ is a $B$-morphism, and whose morphisms from $\begin{pmatrix}
      M_1 \\ 
      M_2   
    \end{pmatrix}_{\varphi^M}$ to $\begin{pmatrix}
      N_1 \\ 
      N_2   
    \end{pmatrix}_{\varphi^N}$ are pairs
    $\begin{pmatrix}
      f_1 \\ 
      f_2   
    \end{pmatrix}$
such that $f_1\in \Hom_A(M_1,N_1)$, $f_2\in
\Hom_B(M_2,N_2)$ satisfying that the following diagram is commutative.

$$ \xymatrix{
    U\otimes_AM_1 \ar[r]^-{\varphi^M} \ar[d]_{1_U\otimes f_1}   &   M_2\ar[d] ^{f_2}
    \\
     U\otimes_AN_1\ar[r]^-{\varphi^N}         &  N_2
\\
}$$ 
Since we have the natural isomorphism $$\Hom_B(U\otimes_A M_1,M_2)\cong \Hom_A(M_1,\Hom_B(U,M_2)),$$there is an alternative way of defining $T$-modules and $T$-homomorphisms in terms of maps  $\widetilde{\varphi^M}:M_1\to \Hom_B(U,M_2)$ given by $\widetilde{\varphi^M}(x)(u)=\varphi^M(u\otimes x)$ for  each $u\in U$ and $x\in M_1.$

Analogously, the category Mod-$T$ of right $T$-modules is equivalent to the category $\Omega_T$ whose objects are triples $M=\begin{pmatrix}
      M_1 , M_2   
    \end{pmatrix}_{\psi^M}$, where $M_1\in$ Mod-$A$, $M_2\in$ Mod-$B$ and $\varphi^M : M_2\otimes_B U\rightarrow M_1$ is an $A$-morphism, and whose morphisms from $\begin{pmatrix}
      M_1  , M_2   
    \end{pmatrix}_{\phi^M}$ to $\begin{pmatrix}
      N_1 , 
      N_2   
    \end{pmatrix}_{\phi^N}$ are pairs
    $\begin{pmatrix}
      f_1 ,
      f_2   
    \end{pmatrix}$
such that $f_1\in \Hom_A(M_1,N_1)$, $f_2\in
\Hom_B(M_2,N_2)$ satisfying that the following diagram

$$ \xymatrix{
    M_2\otimes_B U \ar[r]^-{\phi^M} \ar[d]_{f_2\otimes 1_U}  &   M_1\ar[d] ^{f_1} 
    \\
    M_2\otimes_B U\ar[r]^-{\phi^N}         &  N_1
\\
}$$ is commutative. 

In the rest of the paper we shall identify $T$-Mod (resp. Mod-$T$) with $_T\Omega$ (resp. $\Omega_T$). Consequently, through the paper, a left (resp. right) $T$-module will be a triple $M=\begin{pmatrix}
      M_1 \\ 
      M_2   
    \end{pmatrix}_{\varphi^M}$ (resp. $M=\begin{pmatrix}
      M_1 ,  M_2   
    \end{pmatrix}_{\phi^M}$) and, whenever there is no possible confusion, we shall omit the morphisms $\varphi^M$ and $\phi^M$. For example, $_TT$ is identified with $\begin{pmatrix}
    A\\U\oplus B
    \end{pmatrix}$ and $T_T$ is identified with $\begin{pmatrix}
    A\oplus U, B
    \end{pmatrix}$.
      
    A sequence of left T-modules $0\rightarrow\begin{pmatrix}
M_1\\M_2 
\end{pmatrix}\rightarrow\begin{pmatrix}
M'_1\\M'_2 
\end{pmatrix}\rightarrow\begin{pmatrix}
M''_1\\M''_2 
\end{pmatrix}\rightarrow 0$ is exact if and only if both sequences $0\rightarrow
M_1\rightarrow
M'_1\rightarrow
M''_1\rightarrow 0$ and $0\rightarrow
M_2\rightarrow
M'_2\rightarrow
M''_2\rightarrow 0$ are exact.

Throughout this paper,  $T=\begin{pmatrix}
      A & 0\\ 
      U & B  
    \end{pmatrix}$ will be a generalized triangular matrix ring. Given a $T$-module $M=\begin{pmatrix}
      M_1 \\ 
      M_2   
    \end{pmatrix}_{\varphi^M}$ , the $B$-module $\Coker\varphi^M$ will be denoted as $\overline{M}_2$ and the $A$-module $\Ker \widetilde{\varphi^M}$ as $\underline{M_1}$.  A $T$-module $N=\begin{pmatrix}
N_1\\N_2
    \end{pmatrix}_{\varphi^N}$  is a submodule of $M$ if $N_1$ is a submodule of $M_1$, $N_2$ is a submodule of $M_2$ and $\varphi^M|_{U\otimes_A N_1}=\varphi^N$.

As an interesting and special case of triangular matrix rings, we recall that the $T_2$-extension of a ring $R$ is given by 
$$T(R)=\begin{pmatrix}
R&0\\R&R
\end{pmatrix}$$ and the modules over $T(R)$ are triples $M=\begin{pmatrix}
M_1\\M_2
\end{pmatrix}_{\varphi^M}$ where $M_1$ and $M_2$ are $R$-modules and $\varphi^M:M_1\to M_2$ is an $R$-homomorphism. 

There are some pairs of adjoint functors $(\textbf{p},\textbf{q})$, $(\textbf{q},\textbf{h})$ and $(\textbf{s},\textbf{r})$ between the category $T$-Mod and the product category $A$-Mod $\times$ $B$-Mod which are defined as follows:
\begin{enumerate}
\item $\textbf{p}$\;:\;$A$-Mod $\times$ $B$-Mod$\rightarrow T$-Mod is defined as follows: for each object $(M_1,M_2)$
of $A$-Mod$\times B$-Mod, let $\textbf{p}(M_1,M_2)=\begin{pmatrix}
      M_1 \\ 
      (U\otimes_A M_1)\oplus M_2   
    \end{pmatrix}$ with the obvious map and 
for any morphism $(f_1,f_2)$ in $A$-Mod$\times B$-Mod, let $\textbf{p}(f_1,f_2)=\begin{pmatrix}
      f_1 \\ 
      (1_U\otimes_A f_1)\oplus f_2   
    \end{pmatrix}$.

\item $\textbf{q}$\;:\;$T$-Mod$\rightarrow A$-Mod $\times B$-Mod is defined, for each left $T$-module $\begin{pmatrix}
      M_1 \\ 
       M_2   
    \end{pmatrix}$ as $\textbf{q}(\begin{pmatrix}
      M_1 \\ 
       M_2   
    \end{pmatrix})
=(M_1,M_2)$, and for each morphism $\begin{pmatrix}
      f_1 \\ 
       f_2   
    \end{pmatrix}$ in $T$-Mod as $\textbf{q}(\begin{pmatrix}
      f_1 \\ 
       f_2   
    \end{pmatrix})
=(f_1,f_2)$.

\item $\textbf{h}$\;:\;$A$-Mod $\times$ $B$-Mod$\rightarrow T$-Mod is defined as follows: for each object $(M_1,M_2)$
of $A$-Mod$\times B$-Mod, let $\textbf{h}(M_1,M_2)=\begin{pmatrix}
        M_1 \oplus \Hom_B(U,M_2)\\ 
       M_2 
        
    \end{pmatrix}$ with the obvious map and 
for any morphism $(f_1,f_2)$ in $A$-Mod$\times B$-Mod, let $\textbf{h}(f_1,f_2)=\begin{pmatrix}
       f_1 \oplus \Hom_B(U,f_2)\\ 
       f_2 
        
    \end{pmatrix}$.

\item $\textbf{r}$\;:\;$A$-Mod $\times$ $B$-Mod$\rightarrow T$-Mod is defined as follows: for each object $(M_1,M_2)$
of $A$-Mod$\times B$-Mod, let $\textbf{r}(M_1,M_2)=\begin{pmatrix}
      M_1 \\ 
      M_2   
    \end{pmatrix}$ with the zero map and 
for any morphism $(f_1,f_2)$ in $A$-Mod$\times B$-Mod, let $\textbf{r}(f_1,f_2)=\begin{pmatrix}
      f_1 \\ 
      f_2   
    \end{pmatrix}$.

\item $\textbf{s}$\;:\;$T$-Mod$\rightarrow A$-Mod $\times B$-Mod is defined, for each left $T$-module $\begin{pmatrix}
      M_1 \\ 
       M_2   
    \end{pmatrix}$ as $\textbf{s}(\begin{pmatrix}
      M_1 \\ 
       M_2   
    \end{pmatrix})
=(M_1,\overline{M}_2)$, and for each morphism $\begin{pmatrix}
      f_1 \\ 
       f_2   
    \end{pmatrix}$ in $T$-Mod as $\textbf{s}(\begin{pmatrix}
      f_1 \\ 
       f_2   
    \end{pmatrix})
=(f_1,\overline{f}_2)$, where $\overline{f}_2$ is the induced map.
\end{enumerate}
It is easy to see that \textbf{q} is exact. In particular,  \textbf{p} preserves projective objects and \textbf{h} preserves injective objects. Note that the pairs of adjoint functors $(\textbf{p},\textbf{q})$ and  $(\textbf{q},\textbf{h})$ are defined in \cite{EIT}. In general, the three pairs of adjoint functors defined above can be found in \cite{FGR}.

For a future reference, we list these adjointness isomorphisms:
$$\Hom_T(\begin{pmatrix}
      M_1 \\ 
      (U\otimes_A M_1)\oplus M_2   
    \end{pmatrix},N)\cong \Hom_A(M_1,N_1)\oplus  \Hom_B(M_2,N_2).$$

$$\Hom_T(N,\begin{pmatrix}
      M_1 \\ 
      M_2   
    \end{pmatrix}_0)\cong \Hom_A(N_1,M_1)\oplus  \Hom_B(\overline{N}_2,M_2).$$
    
$$\Hom_T(M,\begin{pmatrix}
      N_1 \oplus \Hom_B(U,N_2)\\ 
       N_2   
    \end{pmatrix})\cong \Hom_A(M_1,N_1)\oplus  \Hom_B(M_2,N_2).$$

Now we recall the characterizations of projective, injective and finitely generated $T$-modules.
\begin{lem}\label{projectives-injectives over T} Let $M=\begin{pmatrix}
M_1\\M_2
\end{pmatrix}_{\varphi^M}$ be a $T$-module.
 \item[(1)](\cite[Theorem 3.1]{HV}) $M$ is projective if and only if $M_1$ is projective in $A$-Mod, $\overline{M}_2=Coker \varphi^M$ is projective in $B$-Mod and $\varphi^M$ is injective.
 \item[(2)](\cite[Proposition 5.1]{HV2} ) $M$ is injective if and only if $\underline{M_1}=\Ker \widetilde{\varphi^M}$ is injective in $A$-Mod, $M_2$ is injective in $B$-Mod and $\widetilde{\varphi^M}$ is surjective.
 \item[(3)] (\cite{GW}) $M$ is finitely generated if and only  if  $M_1$ and $\overline{M}_2$ are finitely generated.
\end{lem}

The following Lemma improves  \cite[Lemma 3.2]{M}.
\begin{lem}\label{Ext} Let $M=\begin{pmatrix}
M_1\\M_2
\end{pmatrix}_{\varphi^M}$ and $N=\begin{pmatrix}
N_1\\N_2
\end{pmatrix}_{\varphi^N}$ be two $T$-modules and $n\geq 1$ be an integer number. Then we have the following natural isomorphisms:
\begin{enumerate}

\item If $\Tor_{1\leq i\leq n}^A(U,M_1)=0$, then $\Ext^n_T(\begin{pmatrix}
M_1\\ U\otimes_AM_1
\end{pmatrix},N)\cong \Ext^n_A(M_1,N_1).$
\item  $\Ext^n_T(\begin{pmatrix}
0\\M_2
\end{pmatrix},N) \cong \Ext^n_B(M_2,N_2).$
\item $\Ext^n_T(M,\begin{pmatrix}
N_1\\0
\end{pmatrix})\cong \Ext^n_A(M_1,N_1).$
\item If $\Ext_B^{1\leq i\leq n}(U,N_2)=0$, then $\Ext^n_T(M,\begin{pmatrix}
\Hom_B(U,N_2)\\N_2
\end{pmatrix})\cong \Ext^n _B(M_2,N_2).$
\end{enumerate}
\end{lem}
\proof We prove only $1.$, since $2.$ is similar and $3.$ and $4.$ are dual. Assume that $\Tor_{1\leq i\leq n}^A(U,M_1)=0$ and consider an exact sequence of $A$-modules 
$$0\to K_1\to P_1\to M_1\to 0$$
where $P_1$ is projective. So, there exists an exact sequence of $T$-modules 
$$0\to \begin{pmatrix}
K_1\\ U\otimes_AK_1
\end{pmatrix}\to \begin{pmatrix}
P_1\\ U\otimes_AP_1
\end{pmatrix}\to \begin{pmatrix}
M_1\\ U\otimes_AM_1
\end{pmatrix}\to 0$$
where $\begin{pmatrix}
P_1\\ U\otimes_AP_1
\end{pmatrix}$ is projective by Lemma \ref{projectives-injectives over T}.

Let $n=1$. By applying the functor $\Hom_T(-,N)$ to the above short exact sequence and since $\begin{pmatrix}
P_1\\ U\otimes_AP_1
\end{pmatrix}$ and $P_1$ are projectives, we get a commutative diagram with exact rows
 
$$\xymatrix{
\Hom_T({\begin{pmatrix}
P_1\\ U\otimes_AP_1
\end{pmatrix}},N) \ar[d]^{\cong}\ar[r]  	&\Hom_T({\begin{pmatrix}
K_1\\ U\otimes_AK_1
\end{pmatrix}},N) \ar[d]^{\cong} \ar@{>>}[r] &\Ext^1_T({\begin{pmatrix}
M_1\\ U\otimes_AM_1
\end{pmatrix}},N) \ar[d]  \\
\Hom_A(P_1,N_1) \ar[r]    &\Hom_A(K_1,N_1) \ar@{>>}[r]  &\Ext^1_A(M_1,N_1) 
}$$

where the first two columns are just the natural isomorphisms given by adjointeness and the last two horizontal rows are epimorphisms. Thus, the induced map $$\Ext^1_T({\begin{pmatrix}
M_1\\ U\otimes_AM_1
\end{pmatrix}},N)\to  \Ext^1_A(M_1,N_1)$$ is an isomorphism such that the above diagram is commutative.
  
Assume that $n>1$. Using the long exact sequence, we get a commutative diagram  with exact rows
$$\xymatrix{
0 \ar[r]  	&\Ext^{n-1}_T({\begin{pmatrix}
K_1\\ U\otimes_AK_1
\end{pmatrix}},N) \ar[d]^{\cong}_{\sigma} \ar[r]_{\cong}^{f} &\Ext^n_T({\begin{pmatrix}
M_1\\ U\otimes_AM_1
\end{pmatrix}},N) \ar@{-->}[d]\ar[r] &0 \\
0 \ar[r]    &\Ext^{n-1}_A(K_1,N_1) \ar[r]^{g}_{\cong}  &\Ext^n_A(M_1,N_1) \ar[r] &0
}$$
where $\sigma$ is a natural isomorphism by induction, since $\Tor_k^A(U,K_1)=0$ for every $k\in\{1,\cdots,n-1\}$ because of the exactness of the following sequence
$$0 = \Tor_{k+1}^A(U,M_1)\to \Tor_k^A(U,K_1)\to \Tor_k^A(U,P_1)= 0.$$
Thus, the composite map $$g\sigma f^{-1}:\Ext^n_T({\begin{pmatrix}
M_1\\ U\otimes_AM_1
\end{pmatrix}},N)\to  \Ext^n_A(M_1,N_1)$$ is a natural  isomorphism, as desired.
\cqfd

Since $T$ can be viewed as a trivial extension (see \cite{FGR} and \cite{BBG} for more details), the following Lemma can be easily deduced from (\cite[Theorem 3.1 and Theorem 3.4]{BBG}). For the convenience of the reader, we give a proof.
\begin{lem}\label{Add-Prod}Let $X=\begin{pmatrix}
X_1\\X_2
\end{pmatrix}_{\varphi^X}$ be a $T$-module and $(C_1,C_2)\in A$-Mod $\times\; B$-Mod.
\begin{enumerate}

   \item $X\in \Add _T(\textbf{p}(C_1,C_2))$ if and only if   \begin{itemize}
   
          \item[(i)]   $X\cong \textbf{p}(X_1,\overline{X}_2)$,
          \item[(ii)]   $X_1\in \Add_A(C_1)$ and $\overline{X}_2 \in \Add_B(C_2).$
   \end{itemize}
   In this case, $\varphi^X$ is injective.
   \item  $X\in \Prod_T(\textbf{h}(C_1,C_2))$ if and only if   \begin{itemize}
   
          \item[(i)]  $X\cong \textbf{h}(\underline{X_1},X_2)$,
          \item[(ii)]   $\underline{X_1}\in \Prod_A(C_1)$ and $X_2 \in \Prod_B(C_2).$
   \end{itemize}
   In this case, $\widetilde{\varphi^X}$ is surjective.
   \end{enumerate}
          
\end{lem}
\proof We only need to prove (1), since (2) is dual.

For the "if" part. If  $X_1\in \Add_A(C_1)$ and $\overline{X}_2 \in \Add_B(C_2)$, then $X_1\oplus Y_1=C^{(I_1)}$ and $\overline{X}_2 \oplus Y_2=C_2^{(I_2)}$ for some $(Y_1,Y_2)\in A$-Mod$\times B$-Mod and some sets $I_1$ and $I_2$.  Without loss of generality, we may assume that $I=I_1=I_2.$ Then 
 
 \begin{eqnarray*}
 X\oplus \textbf{p}(Y_1,Y_2)&\cong & \textbf{p}(X_1,\overline{X}_2)\oplus \textbf{p}(Y_1,Y_2) \\
  &=& \begin{pmatrix}
X_1\\ ( U\otimes_A X_1)\oplus \overline{X}_2 
\end{pmatrix}\oplus \begin{pmatrix}
Y_1\\ ( U\otimes_A Y_1)\oplus Y_2 
\end{pmatrix}\\
&\cong & \begin{pmatrix}
C_1^{(I)}\\ ( U\otimes_A C_1^{(I)})\oplus C^{(I)}_2 
\end{pmatrix}\\
&\cong & \begin{pmatrix}
C_1\\ ( U\otimes_A C_1)\oplus C_2 
\end{pmatrix}^{(I)}\\
&=& \textbf{p}(C_1,C_2)^{(I)}.
 \end{eqnarray*}
Hence  $X\in \Add_T(\textbf{p}(C_1,C_2)).$

Conversely, let  $X\in \Add _T(\textbf{p}(C_1,C_2))$ and $Y=\begin{pmatrix}
Y_1\\Y_2
\end{pmatrix}_{\varphi^Y}$ be a $T$-module such that $X\oplus Y=\textbf{p}(C_1,C_2)^{(I)}$ for some set $I$. Then $\varphi^X$ is injective, as $X$ is a submodule of $C:=\textbf{p}(C_1,C_2)^{(I)}$ and $\varphi^C$ is injective. Consider now the  split exact sequence $$0\to Y\stackrel{\begin{pmatrix}
\lambda_1\\\lambda_2
\end{pmatrix}}{\longrightarrow} C \stackrel{\begin{pmatrix}
p_1\\p_2
\end{pmatrix}}{\longrightarrow} X\to 0 $$
which induces  
the following commutative
diagram with exact rows and colmmuns

$$   \xymatrix{
 &   &  &   & &     \\
0\ar[r] & U\otimes_A Y_1\ar@{^{(}->}[d]_{\varphi^{Y}} \ar[r]^{1_U\otimes\lambda_1}  &  U\otimes_A C_1^{(I)}\ar@{^{(}->}[d]_{\varphi^{C}} \ar[r]^{1_U\otimes p_1}&  U\otimes_A X_1\ar@{^{(}->}[d]_{\varphi^{X}} \ar[r] &  0   \\
 0\ar[r] &Y_2\ar[d]_{\overline{\varphi^{Y}}} \ar[r]^{\lambda_2\hspace{0.9cm}} &  U\otimes_A C_1^{(I)} \oplus C_2^{(I)}\ar[d]_{\overline{\varphi^{C}}} \ar[r]^{\hspace{0.9cm}p_2} &  X_2 \ar[d]_{\overline{\varphi^{X}}}  \ar[r] & 0\\
0\ar[r] & \overline{Y_2}\ar[r]^{\overline{\lambda_2}} \ar[d] &  C_2^{(I)}  \ar[d] \ar[r]^{\overline{p_2}}& \overline{X_2}\ar[r] \ar[d] & 0\\
&0  &0 &0  &    }  $$
where $\overline{\varphi^{X}}$, $\overline{\varphi^{C}}$ and $\overline{\varphi^{X}}$ are the canonical projections.
Clearly, $p_1:C_1^{(I)}\to X_1 $ and $\overline{p_2}:C_2^{(I)}\to \overline{X}_2 $ are split epimorphisms. Then $X_1\in \Add_A(C_1)$ and $\overline{X}_2 \in \Add_B(C_2).$ It remains to prove that $X\cong \textbf{p}(X_1,\overline{X}_2)$.  For this, it suffices to prove that the short exact sequence 
$$0\to U\otimes_A X_1\stackrel{\varphi^X}{\to} X_2\stackrel{\overline{\varphi^X}}{\to} \overline{X}_2 \to 0$$
splits. Let $r_2$ be the retraction of $\overline{p_2}$. If $i:C_2^{(I)}\to (U\otimes_A C_1^{(I)})\oplus C_2^{(I)} $ denotes the canonical injection, then 
$\overline{\varphi^X}p_2ir_2=\overline{p_2}\overline{\varphi^C}ir_2=\overline{p_2}r_2=1_{\overline{X_2}}$ and the proof is finished.
\cqfd

\begin{rem}\begin{enumerate}
\item  Since the class of projective modules over $T$ is nothing but the class $\Add_T(T)$, when we take $C_1=A$ and $C_2=B$ in Lemma \ref{Add-Prod}, we recover the characterization of projective and injective  $T$-modules.
\item Let $(C_1,C_2)\in A$-Mod $\times\; B$-Mod. By Lemma \ref{Add-Prod}, every module in $\Add_T(\textbf{p}(C_1,C_2))$ has the form $\textbf{p}(X_1,X_2)$ for some $X_1\in\Add_A(C_1)$ and  $X_2\in\Add_B(C_2)$
	\end{enumerate}
\end{rem}

\section{w-Tilting modules }
In this section, we study when the functor $\textbf{p}$ preserves w-tilting modules.

It is well known that the functor $\textbf{p}$ preserves projective modules. However,
the functor $\textbf{p}$ does not preserve w-tilting modules in general, as the following example shows.
\begin{exmp}\label{example 1---2} Let $Q$ be the quiver 

$$e_1\stackrel{\alpha}{\longrightarrow} e_2,$$ 
and let $R=kQ$ be the path algebra over an algebraic closed field $k$. Put $P_1=Re_1$, $P_2=Re_2$, $I_1=\Hom_k(e_1R,k)$ and $I_2=\Hom_k(e_2R,k).$ Note that, $C_1$ and $C_2$ are projective and injective $R$-modules, respectively. By \cite[Example 2.3]{ATY},

$$C_1=P_1\oplus P_2(=R)\;\;\;\text{and} \;\;\;C_2=I_1\oplus I_2.$$
are semidualzing $(R,R)$-bimodules and then $w$-tilting $R$-modules. Now, consider the triangular matrix ring $$T(R)=\begin{pmatrix}
R& 0\\R&R
\end{pmatrix}.$$ We claim that  $\textbf{p}(C_1,C_2)$ is not a w-tilting $T(R)$-module. Note that $I_1$ is not projective. Since  $R$ is left hereditary by \cite[Proposition 1.4]{ARS}, $\pd_R(I_1)=1$. Hence $\Ext^1_R(I_1,R)\neq 0.$ Using Lemma $\ref{Ext}$, we get that $\Ext_{T(R)}^1(\textbf{p}(C_1,C_2),\textbf{p}(C_1,C_2))\cong  \Ext_R^1(C_1,C_1)
                \oplus  \Ext_R^1(C_2,C_1)\oplus\Ext_R^1(C_2,C_2)              
                \cong \Ext_R^1(I_1,R) \neq 0               $. Thus, $\textbf{p}(C_1,C_2)$ is not a  w-tilting $T(R)$-module.
\end{exmp} 
Motivated by the definition of compatible bimodules in \cite[Definition 1.1]{Z}, we introduce the following definition which will be crucial throughout the rest of this paper.
\begin{defn}\label{C-compatible}
Let $(C_1,C_2)\in A$-Mod $\times\; B$-Mod and $C=\textbf{p}(C_1,C_2)$. The bimodule  $_BU_A$ is said to be $C$-compatible if the following two conditions hold:
 \begin{enumerate}

\item[(a)]The complex $U\otimes_A \textbf{X}_1$ is exact for every  exact sequence in $A$-Mod
$$\textbf{X}_1:\;\cdots\to P^1_1\to P^0_1\to C^0_1\to C^1_1\to \cdots$$
where the $P^i_1$'s are all projective and $C^i_1\in \Add_A(C_1)$ $\forall i$.

\item[(b)]The complex  $\Hom_B(\textbf{X}_2,U\otimes_A \Add_A(C_1))$ is exact for every\\ $\Hom_B(-,\Add_B(C_2))$-exact exact sequence in $B$-Mod
$$\textbf{X}_2:\;\cdots\to P^1_2\to P^0_2\to C^0_2\to C^1_2\to \cdots$$
where the $P^i_2$'s are all projective and $C^i_2\in \Add_B(C_2)$ $\forall i$.
\end{enumerate}
Moreover, $U$ is called weakly $C$-compatible if it satisfies $(b)$ and the following condition
\begin{enumerate}
\item[(a')]The complex $U\otimes_A \textbf{X}_1$ is exact for every  $\Hom_A(-,\Add_A(C_1))$-exact exact sequence in $A$-Mod
$$\textbf{X}_1:\;\cdots\to P^1_1\to P^0_1\to C^0_1\to C^1_1\to \cdots$$
where the $P^i_1$'s are all projective and $C^i_1\in \Add_A(C_1)$ $\forall i$. 
\end{enumerate}
When $C=\;_TT=\textbf{p}(A,B)$, the bimodule $U$ will be called simply (weakly) compatible. 
\end{defn}
\begin{rem}
 \begin{enumerate}
 \item It is clear by the definition that every $C$-compatible is weakly $C$-compatible.
 \item The  $(B,A)$-bimodule $U$ is weakly compatible if and only if the functor $U\otimes_A-:A$-Mod $\to B$-Mod is weak compatible (see \cite{HZ}). 
 \item  If  $A$ and $B$ are Artin algebras and since $_TT=\begin{pmatrix}
 A\\U\oplus B
 \end{pmatrix}=\textbf{p}(A,B)$, it is easy to see that $_TT$-compatible  bimodules are
 nothing but compatible $(B,A)$-bimodules as defined in \cite{Z}.
 
 \end{enumerate}
\end{rem}

The following  can be applied to produce examples of (weakly) $C$-compatible bimodules later on. 
\begin{lem}Let $C=\textbf{p}(C_1,C_2)$ be a $T$-module.
\begin{enumerate}
\item Assume that $\Tor_1^A(U,C_1)=0$. If $\fd_A(U)<\infty$, then $U$ satisfies  $(a)$.
\item Assume that $\Ext^1_B(C_2,U\otimes_AC^{(I)}_1)=0$ for every set $I$. If $\id_B(U\otimes_AC_1)<\infty$, then $U$ satisfies $(b)$. 

\item If $U\otimes_A C_1\in\Add_B(C_2)$, then $U$ satisfies $(b)$.

\end{enumerate}
\end{lem}
\proof  (3) is clear. We only prove (1), as (2) is similar. Consider an exact sequence of $A$-modules 
$$\textbf{X}_1:\;\cdots\to P^1_1\to P^0_1\to C^0_1\to C^1_1\to \cdots$$
 where the $P^i_1$'s are all projective and $C^i_1\in \Add_A(C_1)$ $\forall i$. We use induction on $\fd_AU$. If $\fd_AU=0,$ then the result is trivial. Now suppose that $\fd_AU=n\geq 1$. Then, there exists an exact sequence of right $A$-modules 
$$0\to L\to F\to U\to 0$$
where $\fd_AL=n-1$ and $F$ is flat. Applying the functor $-\otimes\textbf{X}_1$ to the above short exact sequence, we get the commutative diagram with exact rows

$$   \xymatrix{
:\ar[d] &:\ar[d]  &:\ar[d]  &:\ar[d]  &     \\
0 \ar[d]\ar[r] & L\otimes P^0_1\ar[d] \ar[r]  &  F\otimes_A P^0_1\ar[d] \ar[r]&  U\otimes_A P^0_1\ar[d] \ar[r] &  0   \\
 \Tor^A_1(U,C^0_1)\ar[d]\ar[r]& L\otimes_A C^0_1 \ar[d] \ar[r] & F\otimes_A C^0_1\ar[d] \ar[r] &  U\otimes_A C^0_1 \ar[d]  \ar[r] & 0\\
\Tor^A_1(U,C^1_1) \ar[d]\ar[r]& L\otimes_A C^1_1\ar[r] \ar[d] &  F\otimes_A C^1_1  \ar[d] \ar[r]& U\otimes_A C^1_1\ar[r] \ar[d] & 0\\
 : &: &:  &:    }  $$ 
Since $\Tor_1^A(U,C_1)=0$, the above diagram induces an exact sequence of complexes 
$$0\to L\otimes_A\textbf{X}_1\to F\otimes_A\textbf{X}_1\to U\otimes_A\textbf{X}_1\to 0.$$
By induction hypothesis, the complexes $L\otimes_A\textbf{X}_1$ and $F\otimes_A\textbf{X}_1$ are exact. Thus $U\otimes_A\textbf{X}_1$ is exact as well.   \cqfd

Given a $T$-module $C=\textbf{p}(C_1,C_2)$, we have  simple characterizations of conditions $(a')$ and $(b)$ if $C_1$ and $C_2$ are w-tilting.
\begin{prop}\label{charaterization of condition a and b }
Let $C=\textbf{p}(C_1,C_2)$ be a $T$-module. 
\begin{enumerate}
\item If $C_1$ is w-tilting, then the following assertions are equivalent:
\begin{enumerate}
\item[(i)] $U$ satisfies $(a')$.
\item[(ii)] $\Tor_1^A(U,G_1)=0$, $\forall G_1\in G_{C_1}P(A)$.
\item[(iii)] $\Tor_{i\geq 1}^A(U,G_1)=0$, $\forall G_1\in G_{C_1}P(A)$.

\end{enumerate}In this case, $\Tor_{i\geq 1}^A(U,C_1)=0.$

\item If $C_2$ is w-tilting, then the following assertions are equivalent:
\begin{enumerate}
\item[(i)] $U$ satisfies $(b)$.
\item[(ii)] $\Ext^1_B(G_2,U\otimes_AX_1)=0$, $\forall G_2\in G_{C_2}P(B)$, $\forall X_1\in \Add_A(C_1)$.
\item[(iii)] $\Ext^{i\geq 1}_B(G_2,U\otimes_AX_1)=0$, $\forall G_2\in G_{C_2}P(B)$, $\forall X_1\in \Add_A(C_1)$.
\end{enumerate}In this case, $\Ext^{i\geq 1}_B(C_2,U\otimes_AX_1)=0$, $\forall X_1\in \Add_A(C_1).$
\end{enumerate}
\end{prop}
\proof  We only prove (1), since (2) is similar.

$(i)\Rightarrow (iii)$ Let $G_1\in G_{C_1}P(R)$. There exists a $\Hom_A(-,\Add_A(C_1))$-exact exact sequence in $A$-Mod
$$\textbf{X}_1:\;\cdots\to P^1_1\to P^0_1\to C^0_1\to C^1_1\to \cdots$$
where the $P^i_1$'s are all projective, $G_1\cong \Im(P^0_1\to C^0_1)$ and $C^i_1\in \Add_A(C_1)$ $\forall i$. By condition $(a')$, $U\otimes_A\textbf{X}_1$ is exact, which means in particular that $\Tor^A_{i\geq 1}(U,G_1)=0.$

$(iii)\Rightarrow (ii)$ Clear.

$(ii)\Rightarrow (i)$ Follows by \cite[Corollary 2.13]{BGO1}.

Finally, to prove that  $\Tor_{i\geq 1}^A(U,C_1)=0$, note that  $C_1\in G_{C_1}P(A)$  by \cite[Theorem 2.12]{BGO1}.
\cqfd 

In the following proposition, we study when $\textbf{p}$ preserves w-tilting (tilting) modules.
\begin{prop}\label{when p preserves w-tilting}
Let $C=\textbf{p}(C_1,C_2)$ be a $T$-module and assume that $U$ is weakly $C$-compatible. If $C_1$ and $C_2$ are  w-tilting (tilting), then $\textbf{p}(C_1,C_2)$ is w-tilting (tilting).
\end{prop}
\proof By Lemma \ref{projectives-injectives over T},  the functor $\textbf{p}$ preserves finitely generated modules, so we only need prove the statement for w-tilting.  Assume that $C_1$ and $C_2$ are w-tilting and let $I$ be a set. Then $\Ext_A^{i\geq 1}(C_1,C_1^{(I)})=0$ and $\Ext_B^{i\geq 1}(C_2,C_2^{(I)})=0$. By Proposition above, we have 
$\Ext_B^{i\geq 1}(C_2,U\otimes_AC_1^{(I)})=0$ and $\Tor^A_{i\geq 1}(U,C_1)=0$. Using Lemma \ref{Ext}, for every $n \geq 1$ we get that 
\begin{eqnarray*}
\Ext_T^n(C,C^{(I)})&=&\Ext_T^n(\textbf{p}(C_1,C_2),\textbf{p}(C_1,C_2)^{(I)})\\
&\cong& \Ext_A^n(C_1,C_1^{(I)})\oplus\Ext_B^n(C_2,U\otimes_AC_1^{(I)})\oplus\Ext_B^n(C_2,C_2^{(I)})\\
&=&0.
\end{eqnarray*}
Moreover, there exist exact sequences 
$$\textbf{X}_1:\;0\to A\to C_1^0\to C_1^1\to\cdots$$
and 
$$\textbf{X}_2:\;0\to B\to C_2^0\to C_2^1\to\cdots$$
which are $\Hom_A(-,\Add_A(C_1))$-exact and $\Hom_B(-,\Add_B(C_2))$-exact, respectively, and such that $C_1^i\in \Add_A(C_1)$ and $C_2^i\in \Add_B(C_2)$ for every $i\in\N.$
Since $U$ is weakly $C$-compatible, the complex $U\otimes_A\textbf{X}_1$ is exact. So we construct in $T$-Mod the exact sequence
$$\textbf{p}(\textbf{X}_1,\textbf{X}_2):\;0\to T\to \textbf{p}(C_1^0,C_2^0)\to \textbf{p}(C_1^1,C_2^1)\to\cdots$$
where $\textbf{p}(C_1^i,C_2^i)=\begin{pmatrix}
C_1^i\\(U\otimes_A C_1^i)\oplus C_2^i
\end{pmatrix}\in \Add_T(\textbf{p}(C_1,C_2))$, $\forall i\in\N,$ by Lemma \ref{Add-Prod}(1). 

Let $X\in \Add_T(\textbf{p}(C_1,C_2))$. As a consequence of Lemma \ref{Add-Prod}(1), $X=\textbf{p}(X_1,X_2)$ where $X_1\in \Add_A(C_1)$ and $X_2\in \Add_B(C_2).$ 
Using the adjoitness $(\textbf{p},\textbf{q})$, we get an isomorphism of complexes $$\Hom_T(\textbf{p}(\textbf{X}_1,\textbf{X}_2),X)\cong \Hom_A(\textbf{X}_1,X_1)\oplus \Hom_B(\textbf{X}_2,U\otimes X_1)\oplus \Hom_B(\textbf{X}_2,X_2).$$ But $\Hom_A(\textbf{X}_1,X_1)$ and $\Hom_B(\textbf{X}_2,X_2)$ are  exact and the complex $\Hom_B(\textbf{X}_2,U\otimes X_1)$ is also exact since $U$ is weakly $C$-compatible. Then,  $\Hom_T(\textbf{p}(\textbf{X}_1,\textbf{X}_2),X)$ is exact as well and the proof is finished. \cqfd

Now, we illustrate Proposition \ref{when p preserves w-tilting} with two applications. 
\begin{cor} Let $C=\textbf{p}(C_1,C_2)$ be a $T$-module, $A'$ and $B'$ be two rings  such that  $_AC_{A'}$ and $_AC_{B'}$ are bimodules and assume that $U$ is weakly $C$-compatible. If $_AC_{A'}$ and $_AC_{B'}$ are semidualizing bimodules, then $\textbf{p}(C_1,C_2)$ is a semidualizing $(T,\End_T(C))$-bimodule.
\end{cor}
\proof Follows by Proposition \ref{when p preserves w-tilting} and  \cite[Corollory 3.2]{W}. \cqfd

\begin{cor} \label{application 2}  Let $R$ and $S$ be rings, $\theta:R \rightarrow S$ be a homomorphism with $S_R$ flat, and $T=T(\theta)=:\begin{pmatrix}
	R&0\\S&S
	\end{pmatrix}$. Let $C_1$ be an $R$-module such that $S\otimes_RC_1\in\Add_R(C_1)$ (for instance, if $R$ is commutative or $R=S$). If  $_RC_1$  is w-tilting, then 
\begin{enumerate}
	\item $S\otimes_RC_1$ is a w-tilting $S$-module.
	\item  $C=\begin{pmatrix}
	C_1\\ (S\otimes_RC_1)\oplus (S\otimes_RC_1)
	\end{pmatrix}$ is a w-tilting $T(\theta)$-module.
\end{enumerate}  
\end{cor}
\proof 1. Let $C_2=S\otimes_RC_1$ and note that $C=\textbf{p}(C_1,C_2)$ and that $_SS_R$ is $C$-compatible. So, by Proposition \ref{when p preserves w-tilting}, we only need to prove that $C_2$ is w-tilting $S$-module.

Since $_RC_1$ is w-tilting, there exist $\Hom_R(-,\Add_R(C_1))$-exact exact sequences

$$\textbf{P}:\cdots\to P_1\to P_0\to C\to 0 $$
and 
$$\textbf{X}:0\to R\to C_0\to C_1\to\cdots $$
 with each $_RP_i$ projective and $_RC_i\in\Add_R(C_1)$. Since $S_R$ is flat, we get an exact sequence
$$S\otimes_R\textbf{P}:\cdots\to S\otimes_RP_1\to S\otimes_RP_0\to S\otimes_RC\to 0 $$
and
 $$S\otimes_R\textbf{X}:0\to S\to S\otimes_RC_0\to S\otimes_RC_1\to\cdots $$ 
with each $S\otimes_RP_i$ a projective $S$-module and $S\otimes_RC_i\in\Add_R(C_2)$. 

We prove now that $S\otimes_R\textbf{P}$ and $S\otimes_R\textbf{X}$ are $\Hom_S-(,\Add_S(C_2))$-exact . Let $I$ be a set. Then, the complex
$\Hom_S(S\otimes_R\textbf{P},S\otimes_RC_1^{(I)})\cong  \Hom_R(\textbf{P},\Hom_S(S,S\otimes_RC_1^{(I)}))\cong \Hom_R(\textbf{P},S\otimes_RC_1^{(I)})$ is exact since  $S\otimes_RC_1^{(I)}\in \Add_R(C_1)$. Similarly, $S\otimes_R\textbf{X}$ is $\Hom_S(-,\Add_S(C_2))$-exact.

2. This assertion follows from Proposition \ref{when p preserves w-tilting}. We only need to note that $S$ is weakly $C$-compatible since $S_R$ is flat and $S\otimes_R C_1\in\Add_R(C_2)$.   \cqfd

We end this section with an example of a w-tilting module that is neither projective nor injective. 
\begin{exmp}  Take $R$  and $C_2$ as in example \ref{example 1---2}. So, by Corollary \ref{application 2}, $C=\begin{pmatrix} C_2\\C_2\oplus C_2\end{pmatrix}$ is a w-tilting $T(R)$-module. By Lemma \ref{projectives-injectives over T}, $C$ is not projective since $C_2$ is not and it is not injective since the map $\widetilde{\varphi^C}:C_2\to C_2\oplus C_2$ is not surjective.

Moreover, by \cite[Proposition 2.6]{ARS}, $gl.dim (T(R))=gl.dim(R)+1\leq 2$. So, if $0\to T(R)\to E^0\to E^1\to E^2\to 0$  is an injective resolution of $T(R)$, then $C^1= E^0\oplus E^1\oplus E^2$ is a w-tilting  $T(R)$-module. Note that $T(R)$ has at least three w-tilting modules, $C^1$, $C^2=T(R)$ and $C^3=C$.
\end{exmp} 

\section{Relative Gorenstein projective modules}
In this section, we describe $G_C$-projective modules over $T$. Then we use this description to study when the class of $G_C$-projective $T$-modules is a special precovering class. 

Clearly the functor $\textbf{p}$ preserves projective module. So we start by studying when the functor $\textbf{p}$ also preserves relative Gorenstein projective modules. But first we need the following
\begin{lem}\label{when p preseves particular G_C-projective modules} Let $C=\textbf{p}(C_1,C_2)$ be a $T$-module and $U$ be weakly $C$-compatible. 
\begin{enumerate}
\item If $M_1\in G_{C_1}P(A)$, then $\begin{pmatrix}M_1\\
U\otimes_AM_1
\end{pmatrix}\in G_{C}P(T).$
\item If $M_2\in G_{C_2}P(B)$, then
$\begin{pmatrix}
0\\
M_2
\end{pmatrix}\in G_{C}P(T).$
\end{enumerate}
\end{lem}
\proof  $1.$ Suppose that $M_1\in G_{C_1}P(A)$. There exists a $\Hom_A(-,\Add_A(C_1))$-exact exact sequence 
$$\textbf{X}_1:\;\cdots\to P^1_1\to P_1^0\to C_1^0\to C^1_1\to \cdots$$
where the $P^i_1$'s are all projective, $C^i_1\in \Add_A(C_1)$ $\forall i$ and $M_1\cong \Im(P_1^0\to C_1^0)$. Using the fact that  $U$ is weakly $C$-compatible, we get that the complex $U\otimes_A\textbf{X}_1$ is exact in $B$-Mod, which implies that the complex $\textbf{p}(\textbf{X}_1,0):$
$$
\cdots\to\begin{pmatrix}
P^1_1\\U\otimes_A P^1_1
\end{pmatrix}\to\begin{pmatrix}
P^0_1\\U\otimes_A P^0_1
\end{pmatrix}\to\begin{pmatrix}
C^0_1\\U\otimes_A C^0_1
\end{pmatrix}\to\begin{pmatrix}
C^1_1\\U\otimes_A C^1_1
\end{pmatrix}\to\cdots $$
is exact with $\begin{pmatrix}
M_1\\U\otimes_A M_1
\end{pmatrix}\cong \Im(\begin{pmatrix}
P^0_1\\U\otimes_A P^0_1
\end{pmatrix}\to\begin{pmatrix}
C^0_1\\U\otimes_A C^0_1
\end{pmatrix})$. Clearly, $\textbf{p}(P_1^i,0)=\begin{pmatrix}
P^i_1\\U\otimes_A P^i_1
\end{pmatrix}\in \Proj(T)$ and $\textbf{p}(C_1^i,0)=\begin{pmatrix}
C^i_1\\U\otimes_A C^i_1
\end{pmatrix}\in \Add_T(C)$ $\forall i\in \N$ by Lemma \ref{projectives-injectives over T}(1) and Lemma \ref{Add-Prod}(1). If $X\in \Add_T(C)$, then $X_1\in \Add_A(C_1)$ by Lemma \ref{Add-Prod}(1) and using the adjointness, we get that the complex \\ 
$\Hom_T(\textbf{p}(\textbf{X}_1,0),X)\cong \Hom_A(\textbf{X}_1,X_1)$ is exact. Hence $\begin{pmatrix}
M_1\\U\otimes_A M_1
\end{pmatrix}$ is $G_C$-projective. 

$2.$ Suppose that $M_2$ is $G_{C_2}$-projective. There exists a $\Hom_B(-,\Add_B(C_2))$-exact exact sequence 
$$\textbf{X}_2:\;\cdots\to P^1_2\to P_2^0\to C_2^0\to C^1_2\to \cdots$$
where the $P^i_2$'s are all projective, $C^i_2\in \Add_B(C_2)$ $\forall i$ and $M_2\cong \Im(P_2^0\to C_2^0).$ Clearly the complex 
$$\textbf{p}(0,\textbf{X}_2):
\cdots\to\begin{pmatrix}
0\\ P^1_2
\end{pmatrix}\to\begin{pmatrix}
0\\P^0_2
\end{pmatrix}\to\begin{pmatrix}
0\\ C^0_2
\end{pmatrix}\to\begin{pmatrix}
0\\ C^1_2
\end{pmatrix}\to\cdots $$
is exact with $\begin{pmatrix}
0\\ M_2
\end{pmatrix}\cong \Im(\begin{pmatrix}
0\\ P^1_2
\end{pmatrix}\to\begin{pmatrix}
0\\ C^0_2
\end{pmatrix})$, $\textbf{p}(0,P_2^i)=\begin{pmatrix}
0\\ P^i_2
\end{pmatrix}\in \Proj(T)$ and $\textbf{p}(0,C_2^i)=\begin{pmatrix}
0\\ C^i_2
\end{pmatrix}\in \Add_T(C)$ $\forall i,$ by Lemma \ref{projectives-injectives over T}(1) and Lemma \ref{Add-Prod}(1). Let $X\in \Add_T(C)$. Then, by Lemma \ref{Add-Prod}(1), $X=\textbf{p}(X_1,X_2)$ where  $X_1\in \Add_A(C_1)$ and $X_2\in \Add_B(C_2)$. Using adjointness, we get that 
$$\Hom_T(\textbf{p}(0,\textbf{X}_2),X)\cong \Hom_B(\textbf{X}_2,U\otimes_AX_1)\oplus\Hom_B(\textbf{X}_2,X_2)$$
The complex $\Hom_B(\textbf{X}_2,X_2)$ is exact and since $U$ is weakly $C$-compatible, the complex $\Hom_B(\textbf{X}_2,U\otimes_AX_1)$ is also exact. This means that  $\Hom_T(\textbf{p}(0,\textbf{X}_2),X)$ is exact as well and  $\begin{pmatrix}
0\\M_2
\end{pmatrix}$ is  $G_C$-projective.
\cqfd
\begin{prop}\label{when p preserves G_C-projectives }  Let $C=\textbf{p}(C_1,C_2)$  be a $T$-module. If $U$ is weakly $C$-compatible, then the functor $\textbf{p}$ sends $G_{(C_1,C_2)}$-projectives to $G_C$-projectives. The converse holds provided that $C_1$ and $C_2$ are w-tilting.
	
In particular,  $\textbf{p}$ preserves  Gorenstein projective modules if and only if $U$ is weakly compatible.
\end{prop}
\proof Note that 
$$\textbf{p}(M_1,M_2)=\begin{pmatrix}
M_1\\U\otimes_A M_1
\end{pmatrix}\oplus \begin{pmatrix}
0\\M_2
\end{pmatrix}.$$ So this direction follows from Lemma \ref{when p preseves particular G_C-projective modules} and \cite[Proposition 2.5]{BGO1}. 

Conversely, assume that $C_1$ and $C_2$ are w-tilting. By Proposition \ref{charaterization of condition a and b }, it suffices to prove that $\Tor_1^A(U,G_{C_1}P(A))=0=\Ext_B^1(G_{C_2}P(B),U\otimes_A\Add_A(C_1))$. 

Let $G_1\in G_{C_1}P(A)$. By \cite[Corollary 2.13]{BGO1}, there exits a $\Hom_A(-,\Add_A(C_1))$-exact exact sequence $0\to L_1\stackrel{\imath}{\to} P_1\to G_1\to 0$, where $_AP_1$ is projective and $L_1$ is $G_{C_1}$-projective. Note that $A,C_1\in G_{C_1}P(A)$ and  $B,C_2\in G_{C_2}P(B)$ by Lemma \ref{charac of w-tilitng R-mod}. Then $_TT=\textbf{p}(A,B)$ and  $C=\textbf{p}(C_1,C_2)$ are $G_C$-projective, which imply by Lemma \ref{charac of w-tilitng R-mod}, that $C$ is w-tilting. Moreover  $\begin{pmatrix}
L_1\\U\otimes_AL_1
\end{pmatrix}=\textbf{p}(L_1,0)$ is also $G_C$-projective and by  \cite[Corollary 2.13]{BGO1} there exists a  short exact sequence $$0\to \begin{pmatrix}
L_1\\U\otimes_AL_1
\end{pmatrix}\to \begin{pmatrix}
X_1\\X_2
\end{pmatrix}_{\varphi^X}\to \begin{pmatrix}
H_1\\H_2
\end{pmatrix}_{\varphi^H}\to 0 $$
where $X=\begin{pmatrix}
X_1\\X_2
\end{pmatrix}_{\varphi^X}\in\Add_T(C)$ and $H= \begin{pmatrix}
H_1\\H_2
\end{pmatrix}_{\varphi^H}$ is $G_C$-projective.

Since $X_1\in \Add_A(C_1)$, we have the following commutative diagram with exact rows:

$$ \xymatrix{
0 \ar[r] &  L_1 \ar[r]^{\imath} \ar@{=}[d] & 	P_1 \ar[r] \ar[d]  &   G_1\ar[d] \ar[r] & 0
	\\
0 \ar[r] & L_1\ar[r]	 & X_1\ar[r]         &  H_1 \ar[r] &0
	\\
}$$
So if we apply the functor $U\otimes_A-$ to the above diagram, we get the following commutative diagram with exact rows:

$$ \xymatrix{
	& U\otimes_AL_1 \ar[r]^{1_U\otimes\imath} \ar@{=}[d] & 	U\otimes_AP_1 \ar[r] \ar[d]  &   U\otimes_AG_1\ar[d] \ar[r] & 0
	\\
	& U\otimes_AL_1 \ar[r] \ar@{=}[d]	 & U\otimes_AX_1\ar[r]  \ar[d]       &  U\otimes_AH_1 \ar[r] \ar[d] &0
	\\
	0 \ar[r] & U\otimes_AL_1 \ar[r]  	 & X_2\ar[r]         &  H_2 \ar[r] &0 \\
}$$
The commutativity of this diagram implies that the map $1_U\otimes\imath$ injective, and since $P_1$ is projective, $\Tor_1^A(U,G_1)=0$. 

Now let $G_2\in G_{C_2}P(B)$ and $Y_2\in\Add_A(C_1)$. By hypothesis, $\begin{pmatrix}
0\\ G_2
\end{pmatrix}=\textbf{p}(0,G_2)$ is $G_C$-projective and by Lemma \ref{Add-Prod}, $\begin{pmatrix}
Y_1\\ U\otimes Y_1
\end{pmatrix}=\textbf{p}(Y_1,0)\in\Add_T(C)$. Hence 
$\Ext_B^1(G_2,U\otimes_A Y_1)=\Ext_T^1(\begin{pmatrix}
0\\ G_2
\end{pmatrix},\begin{pmatrix}
Y_1\\ U\otimes Y_1
\end{pmatrix})=0$ by Lemma \ref{Ext} and \cite[Proposition 2.4]{BGO1}.
\cqfd

\begin{thm}\label{structure of $G_C$-projective} Let $M=\begin{pmatrix}
M_1\\M_2
\end{pmatrix}_{\varphi^M}$ and $C=\textbf{p}(C_1,C_2)$ be two $T$-modules. If $U$ is $C$-compatible, then the following assertions are equivalent:
\begin{enumerate}

   \item $M$ is $G_C$-projective.
   \item \begin{itemize}
      \item[(i)] $\varphi^M$ is injective.
      \item[(ii)] $M_1$ is $G_{C_1}$-projective   and $\overline{M}_2:=\Coker\;{\varphi^M}$ is $G_{C_2}$-projective.
      \end{itemize}
\end{enumerate}
In this case, if $C_2$ is $\Sigma$-self-orthogonal, then $U\otimes_A M_1$ is $G_{C_2}$-projective if and only if $M_2$ is $G_{C_2}$-projective.
     
\end{thm}
\proof $2. \Rightarrow 1.$ Since $\varphi^M$ is a injective, there exists an exact sequence in $T$-Mod
$$0\to \begin{pmatrix}
M_1\\U\otimes_A M_1
\end{pmatrix}\to M\to \begin{pmatrix}
0\\\overline{M}_2
\end{pmatrix}\to 0$$
Note that $\begin{pmatrix}
M_1\\U\otimes_A M_1
\end{pmatrix} and  \begin{pmatrix}
0\\\overline{M}_2\end{pmatrix}$ are $G_C$-projective $T$-module by Lemma \ref{when p preseves particular G_C-projective modules}. So, $M$ is $G_C$-projective by \cite[Proposition 2.5]{BGO1}.

$1.\Rightarrow 2.$ There exists a $\Hom_T(-,\Add_T(C))$-exact sequence in $T$-Mod $$\textbf{X}=\cdots\to \begin{pmatrix}
P_1^1\\P_2^1
\end{pmatrix}_{\varphi^{P^1}}\to \begin{pmatrix}
P_1^0\\P_2^0
\end{pmatrix}_{\varphi^{P^0}}\to \begin{pmatrix}
C_1^0\\C_2^0
\end{pmatrix}_{\varphi^{C^0}}\to \begin{pmatrix}
C_1^1\\C_2^1
\end{pmatrix}_{\varphi^{C^1}}\to \cdots$$
where $C^i=\begin{pmatrix}
C_1^i\\C_2^i
\end{pmatrix}_{\varphi^{C^i}}\in \Add_T(C)$, $P^i=\begin{pmatrix}
P_1^i\\P_2^i
\end{pmatrix}_{\varphi^{P^i}}\in \Proj(T)$ $\forall i\in\N$, and such that $M \cong \Im(P^0\to C^0).$ Then, we get the exact sequence 
$$\textbf{X}_1=\cdots\to
P_1^1\to 
P_1^0\to 
C_1^0\to
C_1^1\to \cdots$$
where $C^i_1\in \Add_A(C_1)$, $P^i_1\in \Proj(A)$ $\forall i\in\N$ by Lemma \ref{Add-Prod}(1) and Lemma  \ref{projectives-injectives over T}(1), and such that  $M_1\cong \Im(P^0_1\to C^0_1).$ Since $U$ is $C$-compatible, the complex $U\otimes_A\textbf{X}_1$ is exact with $U\otimes_A M_1\cong Im(U\otimes_AP^0_1\to U\otimes_AC^0_1).$ If $\iota_1:M_1\to C^0_1$ and $\iota_2:M_2\to C^0_2$ are the inclusions, then $1_U\otimes \iota_1$ is injective and the following diagram commutes:
$$   \xymatrix{
 &U\otimes_AM_1\ar[r]^{1_U\otimes \iota_1} \ar[d]_{\varphi^{M}}  &U\otimes_AC^0_1\ar[d]_{\varphi^{C^0}}\\ 
& M_2  \ar[r]^{\iota_2}                 &C^0_2   
   }  $$ 

By Lemma \ref{Add-Prod}(1), $\varphi^{C^0}$ is injective, then $\varphi^M$ is also injective.
For every $i \in\N$, $\varphi^{P^i}$ and $\varphi^{C^i}$  are injective by Lemma \ref{projectives-injectives over T} and Lemma \ref{Add-Prod}(1). Then the following diagram with exact
columns

$$   \xymatrix{
 &0\ar[d]  &0\ar[d]  &0\ar[d]  &0\ar[d]   &  &     \\
\cdots \ar[r]& U\otimes_AP_1^1\ar[r] \ar[d]_{\varphi^{P^1}} & U\otimes_AP_1^0\ar[d]_{\varphi^{P^0}} \ar[r]  &  U\otimes_A C_1^0\ar[d]_{\varphi^{C^0}} \ar[r]&  U\otimes_A C_1^1\ar[d]_{\varphi^{C^1}} \ar[r] &  \cdots   \\
\cdots \ar[r]& P_2^1\ar[r] \ar[d]& P_2^0\ar[d] \ar[r] & C_2^0\ar[d] \ar[r] &  C^1_2 \ar[d]  \ar[r] & \cdots\\
\cdots\ar[r] & \overline{P^1_2} \ar[r] \ar[d]& \overline{P^0_2}\ar[r] \ar[d] &  \overline{C^0_2}  \ar[d] \ar[r]& \overline{C_2^1}\ar[r] \ar[d] & \cdots\\
&0 &0  &0 &0  &    }  $$
 is commutative. Since the first row and the second row are exact, we get the exact sequence of B-modules
$$\overline{\textbf{X}}_2:\;\cdots\to \overline{P^1_2} \to \overline{P^0_2}\to  \overline{C^0_2}  \to \overline{C_2^1}\to  \cdots$$
where $\overline{P^i_2}\in \Proj(B)$, $\overline{C^i_2}\in \Add_B(C_2)$ by Lemma \ref{projectives-injectives over T} and Lemma \ref{Add-Prod}(1), and such that $\overline{M}_2=\Im(\overline{P^0_2}\to  \overline{C^0_2})$. It remains to see that $\textbf{X}_1$ and $\overline{\textbf{X}}_2$ are $\Hom_A(-,\Add(C_1))$-exact and $\Hom_B(-,\Add_B(C_2))$-exact, respectively. Let $X_1\in \Add_A(C_1)$ and $X_2\in \Add_B(C_2)$.  Then $\textbf{p}(X_1,0)=\begin{pmatrix}
X_1\\U\otimes_A X_1
\end{pmatrix}\in \Add_T(C)$ and $\textbf{p}(0,X_2)=\begin{pmatrix}
0\\X_2
\end{pmatrix}\in \Add_T(C)$ by Lemma \ref{Add-Prod}(1).  So, by using adjointness, we get that $\Hom_B(\overline{\textbf{X}}_2,X_2)\cong \Hom_T(\textbf{X},\begin{pmatrix}
0\\X_2
\end{pmatrix})$ is exact. Using  adjointness again we get that  
$$\Hom_T(\textbf{X},\begin{pmatrix}
0\\U\otimes_A X_1
\end{pmatrix})\cong \Hom_B(\overline{\textbf{X}}_2,U\otimes_A X_1)$$
and $$\Hom_T(\textbf{X},\begin{pmatrix}
X_1\\0
\end{pmatrix})\cong \Hom_A(\textbf{X}_1, X_1).$$
Note that $C^i\cong \textbf{p}(C^i_1,\overline{C^i_2})$ by Lemma \ref{Add-Prod}(1). Hence $\Ext_T^1(C^i,\begin{pmatrix}
0\\U\otimes_A X_1
\end{pmatrix})\cong \Ext_B^1(\overline{C^i_2},U\otimes_A X_1)=0$ by Lemma \ref{Ext}. So, if we apply the functor $\Hom_T(\textbf{X},-)$ to the sequence 
$$0\to \begin{pmatrix}
0\\U\otimes_A X_1
\end{pmatrix}\to \begin{pmatrix}
X_1\\U\otimes_A X_1
\end{pmatrix} \to \begin{pmatrix}
X_1\\0
\end{pmatrix}\to 0,$$
we get the following exact sequence of complexes

$$0\to\Hom_B(\overline{\textbf{X}}_2,U\otimes_A X_1)   \to  \Hom_T(\textbf{X},\begin{pmatrix}
X_1\\U\otimes_A X_1
\end{pmatrix})  \to   \Hom_A(\textbf{X}_1, X_1)\to   0.$$
Since $U$ is $C$-compatible, it follows that $\Hom_B(\overline{\textbf{X}}_2,U\otimes_A X_1)$ is exact and since $C$ is w-tilting, $ \Hom_T(\textbf{X},\begin{pmatrix}
X_1\\U\otimes_A X_1
\end{pmatrix})$is also exact. Thus $\Hom_A(\textbf{X}_1, X_1)$ is exact and the proof is finished.\cqfd

The following consequence of the above theorem gives the converse of Proposition \ref{when p preserves w-tilting}.
\begin{cor} \label{cns for p(C_1,C_2) to be w-tilting}Let $C=\textbf{p}(C_1,C_2)$ and assume that $U$ is $C$-compatible.  Then $C$ is w-tilting if and only if $C_1$ and $C_2$ are w-tilting.
	
\proof An easy application of  Proposition \ref{charac of w-tilitng R-mod} and Theorem \ref{structure of $G_C$-projective} on the $T$-modules $C=\begin{pmatrix}
C_1\\(U\otimes_A C_1)\oplus C_2
\end{pmatrix}$ and $_TT=\begin{pmatrix}
A\\U\oplus B
\end{pmatrix}$.
\cqfd
One would like to know if every w-tilting $T$-module has the form $\textbf{p}(C_1,C_2)$ where $C_1$ and $C_2$ are w-tilting. The following example gives a negative answer to this question.
\begin{exmp}Let $R$ be a quasi-Frobenius ring and $T(R)=\begin{pmatrix}
R&0\\R&R
	\end{pmatrix}$. Consider the exact sequence of $T$-modules 
	$$0\to T\to \begin{pmatrix}
	R\oplus R\\ R\oplus R
	\end{pmatrix}\to \begin{pmatrix}
	R\\ 0
	\end{pmatrix}\to 0.$$
By Lemma \ref{projectives-injectives over T},  $I^0=\begin{pmatrix}
R\oplus R\\ R\oplus R
\end{pmatrix}$ and $I^1=\begin{pmatrix}
R\\ 0
\end{pmatrix}$ are both injective $T(R)$-modules. Note that $T(R)$ is noetherian (\cite[Proposition 1.7]{GW}) and then we can see that $C:=I^0\oplus I^1$ is a w-tilting $T(R)$-module but does not have the form $\textbf{p}(C_1,C_2)$ where $C_1$ and $C_2$ are w-tilting  by Lemma \ref{Add-Prod} since $I^1\in\Add_{T(R)}(C)$ and $\varphi^{I^1}$ is not injective.

\end{exmp}

\end{cor}
As an immediate consequence of Theorem \ref{structure of $G_C$-projective}, we have the following.

\begin{cor}  Let $R$ be a ring and $T(R)=\begin{pmatrix}
R&0\\R&R
\end{pmatrix}$. If $M=\begin{pmatrix}
M_1\\M_2
\end{pmatrix}_{\varphi^M}$ and $C=\textbf{p}(C_1,C_1)$ are two $T(R)$-modules with $C_1$ $\Sigma$-self-orthogonal, then the following assertions are equivalent: 
\begin{enumerate}
\item $M$ is $G_C$-projective $T(R)$-module
\item $M_1$ and $\overline{M}_2$ are $G_{C_1}$-projective $R$-modules and $\varphi^M$ is injective
\item $M_1$ and $M_2$ are $G_{C_1}$-projective $R$-modules and $\varphi^M$ is injective
\end{enumerate}
\end{cor}

An Artin algebra $\Lambda$ is called Cohen-Macaulay free, or simply, CM-free if any finitely generated Gorenstein projective module is projective. The authors in \cite{EIT}, extended this definition to arbitrary rings and defined strongly  CM-free as rings over which every Gorenstein projective module is projective.
Now, we introduce a relative notion of these rings and give a characterization of when $T$ is such rings.

\begin{defn} Let $R$ be a ring. Given an $R$-module $C$, $R$ is called CM-free (relative to $C$) if $G_CP(R)\cap R\text{-mod}=\add_R(C)$ and it is called strongly CM-free (relative to $C$) if $G_CP(R)=\Add_R(C)$.
\end{defn}
\begin{rem} Let $R$ be a ring and $C$ a $\Sigma$-self-orthogonal $R$-module. Then $\Add_R(C)\subseteq G_CP(R)$ and $\add_R(C)\subseteq G_CP(R)\cap R\text{-mod}$ by \cite[Proposition 2,5, 2,6 and Corollary 2.10]{BGO1}, then  $R$ is CM-free (relative to $C$) if and only if every finitely generated $G_C$-projective is in $\add_R(C)$ and it is strongly CM-free (relative to $C$) if every $G_C$-projective is in $\Add_R(C).$
\end{rem}
Using the above results we refine and extend \cite[Theorem 4.1]{EIT} to our setting.  Note that the condition $B$ is left Gorenstein regular is not needed.
\begin{prop}\label{when T is relative CM-free} Let $_AC_1$ and $_BC_2$ be $\Sigma$-self-orthogonal, and  $C=\textbf{p}(C_1,C_2)$. Assume that $U$ is weakly $C$-compatible and consider the following assertions:
	\begin{enumerate}
	\item $T$ is (strongly) CM-free relative to $C$.
	\item $A$ and $B$ are (strongly) CM-free relative to $C_1$ and $C_2$, respectively.
\end{enumerate}
Then  $1.\Rightarrow 2.$.  If  $U$ is $C$-compatible, then $1.\Leftrightarrow 2.$
\end{prop}
\proof We only prove the the result for relative strongly CM-free, since the the case of relative CM-free is similar. 

$1.\Rightarrow 2.$ By the remark above, we only need to prove that $G_{C_1}P(A)\subseteq \Add_A({C_1})$ and $G_{C_2}P(B)\subseteq \Add_B({C_2})$. Let $M_1$ be a $G_{C_1}$-projective $A$-module and $_BM_2$ a $G_{C_2}$-projective $B$-module. By the assumption and Proposition \ref{when p preserves G_C-projectives }, $\textbf{p}(M_1,M_2)\in G_CP(T)=\Add_T(C)$. Hence $M_1\in\Add_A(C_1)$ and $M_2\in\Add_B(C_2)$ by Lemma \ref{Add-Prod}.

$2.\Rightarrow 1.$ Assume  $U$ is $C$-compatible. Clearly, $C$ is $\Sigma$-self-orthogonal, then by Remark above, we only need to prove that $G_CP(T)\subseteq \Add_T(C)$. Let $M=\begin{pmatrix}
M_1\\M_2
\end{pmatrix}_{\varphi^M}$ be a $G_C$-projective $T$-module. By the assumption and Theorem \ref{structure of $G_C$-projective}, $M_1\in G_{C_1}P(A)=\Add_A(C_1)$ and $\overline{M}_2\in G_{C_2}P(B)=\Add_B(C_2)$ and the map $\varphi^M$ is injective. By the assumption, we can easily see that $\Ext_B^{i\geq 1}(U\otimes_A M_1,\overline{M}_2)=0$. So the map $0\to U\otimes_A M_1\stackrel{\varphi^M }{\to}M_2\to \overline{M}_2\to 0$ splits. Hence $M\cong \textbf{p}(M_1,\overline{M}_2)\in\Add_T(C)$ by Lemma \ref{Add-Prod}.
 \cqfd

Our aim now is to study  special $G_CP(T)$-precovers in $T$-Mod. We start with the following result.
\begin{prop} \label{when a module has precover} Let $C=\textbf{p}(C_1,C_2)$ be w-tilting, $U$ $C$-compatible,  $M=\begin{pmatrix}
M_1\\M_2
\end{pmatrix}_{\varphi^M}$ and  $G=\begin{pmatrix}
G_1\\G_2
\end {pmatrix}_{\varphi^G}$ two $T$-modules with $G$ $G_C$-projective. Then $$f=\begin{pmatrix}
f_1\\f_2
\end{pmatrix}:G\longrightarrow M$$ is a special $G_CP(T)$-precover if and only if 
\begin{enumerate}

\item[(i)] $G_1\stackrel{f_1}{\to} M_1$ is a special $G_{C_1}P(A)$-precover.
\item [(ii)] $G_2\stackrel{f_2}{\to} M_2$ is surjective with its kernel lies 
in $ G_{C_2}P(B)^{\perp_1}$.
\end{enumerate}
In this case, if $G_2\in G_{C_2}P(B)$, then $G_2\stackrel{f_2}{\to} M_2$ is a special $G_{C_2}P(B)$-precover.
\end{prop}
\proof First of all, let $K=\Ker f=\begin{pmatrix}
K_1\\K_2
\end{pmatrix}_{\varphi^K}$ and note that, since $C_1$ is w-tilting, $\Tor_1^A(U,H_1)=0$ for every $H_1\in G_{C_1}P(A)$ by Proposition \ref{charaterization of condition a and b }(1).

$ \Rightarrow $ Since $f$ is surjective, so are $f_1$ and $f_2$. Let $H_1\in G_{C_1}P(A)$ and $H_2\in G_{C_2}P(B)$. Then $\begin{pmatrix}
H_1\\U\otimes_AH_1
\end{pmatrix},\begin{pmatrix}
0\\H_2
\end{pmatrix}\in G_CP(T)$ by Theorem \ref{structure of $G_C$-projective}. Using Lemma \ref{Ext} and the fact that $K$ lies in $G_CP(R)^{\perp_1}$ , we get that $$\Ext^1_A(H_1,K_1)\cong \Ext^1_T(\begin{pmatrix}
H_1\\U\otimes_AH_1
\end{pmatrix},K)=0$$ and  $$\Ext^1_B(H_2,K_2)\cong \Ext^1_T(\begin{pmatrix}
0\\H_2
\end{pmatrix},K)=0.$$ It remains to see that $G_1\in G_{C_1}P(A)$, which is true by Theorem \ref{structure of $G_C$-projective}, since $G$ is $G_C$-projective.

$\Leftarrow $ The morphism $f$ is surjective since $f_1$ and $f_2$ are. So we only need to prove that $K$ lies in $G_CP(R)^{\perp_1}$. Let $H\in G_CP(R)$. By Theorem \ref{structure of $G_C$-projective}, we have the short exact sequence of $T$-modules
$$0\to\begin{pmatrix}
H_1\\ U\otimes_A H_1
\end{pmatrix}\to H\to \begin{pmatrix}
0\\\overline{H}_2
\end{pmatrix}\to 0$$
where $H_1$ is $G_{C_1}$-projective and $\overline{H}_2$ is $G_{C_2}$-projective. So by hypothesis and Lemma \ref{Ext}, we get that 
$\Ext^1_T(\begin{pmatrix}
H_1\\U\otimes_AH_1
\end{pmatrix},K)\cong \Ext^1_A(H_1,K_1)=0$ and  $\Ext^1_T(\begin{pmatrix}
0\\\overline{H}_2
\end{pmatrix},K)\cong \Ext^1_B(\overline{H}_2,K_2) =0$. Then, the exactness of this sequence 
$$\Ext^1_T(\begin{pmatrix}
H_1\\U\otimes_AH_1
\end{pmatrix},K)\to \Ext^1_T(H,K)\to \Ext^1_T(\begin{pmatrix}
0\\\overline{H}_2
\end{pmatrix},K)$$
implies that $\Ext^1_T(H,K)=0.$\cqfd
\begin{thm}\label{the class of G_C-proj is a special precovering}   Let $C=\textbf{p}(C_1,C_2)$ be w-tilting and $U$ $C$-compatible. Then the class $G_CP(T)$ is special precovering in $T$-Mod if and only if  the classes $G_{C_1}P(A)$ and $G_{C_2}P(B)$ are special precovering in $A$-Mod and $B$-Mod, respectively.
\end{thm}
\proof $\Rightarrow $ Let $M_1$ be an $A$-module and $\begin{pmatrix}
G_1\\G_2
\end{pmatrix}_{\varphi^G}\rightarrow \begin{pmatrix}
M_1\\0
\end{pmatrix}$ be a special $G_CP(T)$-precover in $T$-Mod. Then by Proposition \ref{when a module has precover}, $G_1\to M_1$ is a special $G_{C_1}P(A)$-precover in $A$-Mod. 

Let $M_2$ be a $B$-module and $\begin{pmatrix}
0\\ f_2
\end{pmatrix}:\begin{pmatrix}
G_1\\G_2
\end{pmatrix}_{\varphi^G}\rightarrow \begin{pmatrix}
0\\M_2
\end{pmatrix}$ be a special $G_CP(T)$-precover in $T$-Mod. By Proposition \ref{when a module has precover},  $G_1\to 0$ is a special $G_{C_1}P(A)$-precover. Then 
$\Ext^1_A(G_{C_1}P(A),G_1)=0$. On the other hand, by \cite[Proposition 2.8]{BGO1}, there exists an exact sequence of $A$-modules 
$$0\to G_1\to X_1\to H_1\to 0$$
where $X_1\in \Add_A(C_1)$ and $H_1$ is $G_{C_1}$-projective. But this sequence splits, since $\Ext^1_A(H_1,G_1)=0$, which implies that $G_1\in \Add_A(C_1)$. Let $K=\begin{pmatrix}
K_1\\K_2
\end{pmatrix}_{\varphi^K}$ be the kernel of $\begin{pmatrix}
0\\ f_2
\end{pmatrix}$. Note that $K_1=G_1$. So, there exists a commutative diagram

$$   \xymatrix{
 &  &  &   & &     \\
0\ar[r] & U\otimes_A G_1\ar[d]_{\varphi^{K}} \ar@{=}[r]  &  U\otimes_A G_1\ar[d]_{\varphi^G} \ar[r]&  0 \ar[d] \ar[r] &  0   \\
 0\ar[r]& K_2\ar[d] \ar[r] & G_2\ar[d] \ar[r] &  M_2 \ar@{=}[d]  \ar[r] & 0\\
& \overline{K}_2\ar[r] \ar[d] &  \overline{G}_2  \ar[d] \ar[r]& M_2\ar[r] \ar[d] & 0\\
&0  &0 &0  &    }  $$
Using the snake lemma, there exists an exact sequence of $B$-modules 
$$0\to \overline{K}_2\to  \overline{G}_2  \to M_2 \to 0 $$
where $\overline{G}_2$ is $G_{C_2}$-projective by Theorem \ref{structure of $G_C$-projective}. It remains to see that $\overline{K}_2$ lies in $G_{C_2}P(B)^{\perp_1}$. Let $H_2\in G_{C_2}P(B)$. Then 
$\Ext_B^1(H_2,K_2)=0$ by Proposition \ref{when a module has precover} and  $\Ext_B^{i\geq 1}(H_2,U\otimes_A G_1)=0$ by Proposition \ref{charaterization of condition a and b }(2).  From the above diagram, $\varphi^K$ is injective. So, if we apply the functor $\Hom_B(H_2,-)$ to the short exact sequence
$$0\to U\otimes_A G_1\to K_2\to \overline{K}_2\to 0,$$
we get an exact sequence $$\Ext_B^1(H_2,K_2)\to \Ext_B^1(H_2,\overline{K}_2)\to \Ext_B^2(H_2,U\otimes_A G_1)$$
which implies that $\Ext_B^1(H_2,\overline{K}_2)=0$. 

$\Leftarrow$ Note that the functor $U\otimes_A-:A$-Mod$\to B$-Mod is $G_{C_1}P(A)$-exact since $\Tor_1^A(U,G_{C_1}P(A))=0$ by Proposition \ref{charaterization of condition a and b }. So the this direction follows by \cite[Theorem 1.1.]{HZ} since $G_CP(T)=\{M=\begin{pmatrix}
M_1\\M_2
\end{pmatrix}_{\varphi^M}\in T\text{-Mod}|M_1\in G_{C_1}P(A), \;\overline{M_2}\in G_{C_2}P(B) \text{ and $\varphi^M$ is injective  }\}$ by  Theorem \ref{structure of $G_C$-projective}. \cqfd
\begin{cor} Let $R$ be a ring, $T(R)=\begin{pmatrix}
R&0\\R&R
\end{pmatrix}$ and $C=\textbf{p}(C_1,C_1)$ a w-tilting $T(R)$-module. Then $G_CP(T(R))$ is a special precovering class if and only if $G_{C_1}P(R)$ is a special precovering class.
\end{cor}
\section{Relative global Gorenstein dimension}
In this section, we investigate $G_C$-projective dimension of $T$-modules and the left $G_C$-projective global dimension of $T$.

Let $R$ be a ring. Recall (\cite{BGO1}) that a module $M$ is said to have $G_C$-projective dimension less than or equal to $n$, $\Gpc(M)\leq n$, if there is an exact sequence $$0\to G_n\to\cdots\to G_0\to M\to 0$$
with $G_i\in G_CP(R)$ for every $i\in \{0,\cdots,n\}$. If $n$ is the least nonnegative integer for
which such a sequence exists then $\Gpc(M)=n$, and if there is no such $n$ then $\Gpc(M)=\infty$. 

The left $G_C$-projective global dimension of $R$ is defined as: 
$$G_C-PD(R)=sup\{\Gpc(M) \;|\; \text{$M$ is an $R$-module}\}$$ 

\begin{lem} \label{dimension of special modules} Let $C=\textbf{p}(C_1,C_2)$ be w-tilting and $U$ $C$-compatible.
\begin{enumerate}
\item $\2Gpc(M_2)=\Gpc(\begin{pmatrix}
0\\M_2
\end{pmatrix}).$
\item $\1Gpc(M_1)\leq\Gpc(\begin{pmatrix}
M_1\\U\otimes_A M_1
\end{pmatrix}),$ and the equality holds if\\$\Tor^A_{i\geq 1}(U,M_1)=0.$
\end{enumerate}
\end{lem}
\proof 1. Let $n\in \N$ and consider an exact sequence of $B$-modules
$$0\to K_2^n\to G_2^{n-1}\to\cdots\to G_2^0\to M_2 \to 0$$
 where each $G_2^i$ is $G_{C_2}$-projective. Thus, there exists an exact sequence of $T$-modules 
$$0\to \begin{pmatrix}
0\\K_2^n
\end{pmatrix}\to \begin{pmatrix}
0\\G_2^{n-1}
\end{pmatrix}\to\cdots\to \begin{pmatrix}
0\\G_2^0
\end{pmatrix}\to \begin{pmatrix}
0\\M_2
\end{pmatrix}\to 0$$
where each $\begin{pmatrix}
0\\G_2^i
\end{pmatrix}$ is $G_C$-projective by Theorem \ref{structure of $G_C$-projective}. Again, by Theorem \ref{structure of $G_C$-projective},  $\begin{pmatrix}
0\\K_2^n
\end{pmatrix}$ is $G_C$-projective if and only if $K_2^n$ is $G_{C_1}$-projective which means that  $\Gpc(\begin{pmatrix}
0\\M_2
\end{pmatrix})\leq n$ if and only if $\2Gpc(M_2)\leq n$ by \cite[Theorem 3.8]{BGO1}. Hence $\Gpc(\begin{pmatrix}
0\\M_2
\end{pmatrix})=\2Gpc(M_2).$

2. We may assume that $n=\Gpc(\begin{pmatrix}
M_1\\U\otimes_A M_1
\end{pmatrix})<\infty$. By Definition, there exists an exact sequence of $T$-modules
$$0\to G^n\to G^{n-1}\to\cdots\to G^0\to \begin{pmatrix}
M_1\\U\otimes_A M_1
\end{pmatrix}\to 0$$
where each $G^i=\begin{pmatrix}
G_1^i\\G_2^i
\end{pmatrix}_{\varphi^{G^i}}$ is $G_C$-projective. Thus, there exists an exact sequence of $A$-modules 
$$0\to G_1^n\to G_1^{n-1}\to\cdots\to G_1^0\to M_1\to 0$$
where each $G_1^i$ is $G_{C_1}$-projective by Theorem \ref{structure of $G_C$-projective}. So, $\1Gpc(M_1)\leq n$.
Conversely, we prove that $\Gpc(\begin{pmatrix}
M_1\\U\otimes_A M_1
\end{pmatrix})\leq \1Gpc(M_1)$. We may assume that $m:=\1Gpc(M_1)<\infty.$ The hypothesis means that if 
$$\textbf{X}_1 \;: \;0\to K_1^m\to P_1^{m-1}\to\cdots\to P_1^0\to M_1\to 0$$
is an exact sequence of $A$-modules where each $P_1^i$ is projective, then the complex $U\otimes_A\textbf{X}_1$ is exact. Since $C_1$ is w-tilting, each $P_i$ is $G_{C_1}$-projective by \cite[Proposition 2.11]{BGO1} and then $K^m$ is $G_{C_1}$-projective by \cite[Theorem 3.8]{BGO1}. Thus, there exists and exact sequence of $T$-modules 
$$0\to \begin{pmatrix}
K_1^m\\ U\otimes_AK_1^m
\end{pmatrix}\to\begin{pmatrix}
P_1^{m-1}\\ U\otimes_AP_1^{m-1}
\end{pmatrix}\to\cdots \to\begin{pmatrix}
P_1^0\\ U\otimes_AP_1^0
\end{pmatrix}\to\begin{pmatrix}
M_1\\ U\otimes_AM_1
\end{pmatrix}\to 0$$
where $\begin{pmatrix}
K_1^m\\ U\otimes_AK_1^m
\end{pmatrix}$ and all $\begin{pmatrix}
P_1^i\\ U\otimes_AP_1^i
\end{pmatrix}$ are $G_C$-projectives by Theorem \ref{structure of $G_C$-projective}. Therefore, $\Gpc(\begin{pmatrix}
M_1\\U\otimes_A M_1
\end{pmatrix})\leq m=\1Gpc(M_1)$. 
 \cqfd
Given a $T$-module $C=\textbf{p}(C_1,C_2)$, we introduce a strong notion of $G_{C_2}$-projective global dimension of $B$, which will be crucial when we estimate the $G_C$-projective of a $T$-module and the left global $G_C$-projective dimension of $T$. Set $$SG_{C_2}-PD(B)=sup\{\2Gpc_B(U\otimes_A G)\;|\; G\in G_{C_1}P(A)\}.$$

\begin{rem}
\begin{enumerate}
\item Clearly,   $ SG_{C_2}-PD(B)\leq G_{C_2}-PD(B)$.

\item Note that $\pd_B(U)=sup\{\pd_B(U\otimes_A P)\;|\; _AP  \text{ is projective }\; \}$. Therefore, in the classical case, the strong left global dimension of $B$ is nothing but the projective dimension of $_BU$.
\end{enumerate}
\end{rem}
 
\begin{thm}\label{Estimation of dimension of a T-module}
Let $C=\textbf{p}(C_1,C_2)$ be w-tilting, $U$ $C$-compatible, $M=\begin{pmatrix}
M_1\\M_2
\end{pmatrix}_{\varphi^M}$ a $T$-module and $SG_{C_2}-PD(B)<\infty$ . Then 
 
 $$max\{\1Gpc_A(M_1),(\2Gpc_B(M_2))-(SG_{C_2}-PD(B))\}$$
 $$\leq \Gpc(M)\leq$$
 $$ max\{(\1Gpc_A(M_1))+(SG_{C_2}-PD(B))+1,\2Gpc_B(M_2)\}$$
\end{thm}
\proof First of all, note that $C_1$ and $C_2$ are w-tilting by Proposition \ref{cns for p(C_1,C_2) to be w-tilting} and let $k:=SG_{C_2}-PD(B).$ 

Let us first prove that $$max\{\1Gpc(M_1),\2Gpc(M_2)-k\}\leq \Gpc(M).$$ We may assume that $n:=\Gpc(M)<\infty.$  Then, there exists an exact sequence of $T$-modules
$$0\to G^n\to G^{n-1}\to\cdots\to G^0\to M\to 0$$
where each $G^i=\begin{pmatrix}
G_1^i\\G_2^i
\end{pmatrix}_{\varphi^{G^i}}$ is $G_C$-projective.
Thus, there exists an exact sequence of $A$-modules 
$$0\to G_1^n\to G_1^{n-1}\to\cdots\to G_1^0\to M_1\to 0$$
where each $G_1^i$ is $G_{C_1}$-projective by Theorem \ref{structure of $G_C$-projective}. So, $\1Gpc(M_1)\leq n$. By Theorem \ref{structure of $G_C$-projective}, for each $i$, there exists an exact sequence of $B$-modules
 $$0\to U\otimes_AG_1^i\to G_2^i\to \overline{G_2^i}\to 0$$
where $\overline{G_2^i}$ is $G_{C_2}$-projective. Then $\2Gpc(G^i_2)=\2Gpc(U\otimes_AG_1^i)\leq k$ by \cite[Proposition 3.11]{BGO1}. So, using the exact sequence of $B$-modules
$$0\to G_2^n\to G_2^{n-1}\to\cdots\to G_2^0\to M_2\to 0$$ and \cite[Proposition 3.11(4)]{BGO1}, we get that 
$\2Gpc(M_2)\leq n+k.$

Next we prove that  $$\Gpc(M)\leq max\{\1Gpc(M_1)+k+1,\2Gpc(M_2)\}$$ 
We may assume that $$m:=max\{\1Gpc(M_1)+k+1,\2Gpc(M_2)\}< \infty.$$ Then 
$n_1:=\1Gpc(M_1)<\infty$ and $ n_2:=\2Gpc(M_2)<\infty$. Since $\1Gpc(M_1)$\\$=n_1\leq m-k-1$, there exists an exact sequence of $A$-modules 
$$0\to G_1^{m-k-1}\to \cdots \to G_1^{n_2-k}\to\cdots\stackrel{f_1^1}{\to} G_1^0\stackrel{f_1^0}{\to} M_1\to 0$$
where each $G_1^i$ is $G_{C_1}$-projective. Since $C_2$ is w-tilting, there exists an exact sequence of $B$-modules $G_2^0\stackrel{g_2^0}{\to} M_2\to 0$ where $G_2^0$ is $G_{C_2}$-projective by \cite[Corollary 2.14]{BGO1}. Let $K_1^i=\Ker f_1^i$ and define the map $f^0_2:U\otimes_A G_1^0\oplus G_2^0\to M_2$ to be $(\varphi^M(1_U\otimes f^0_1))\oplus g_2^0$. Then, we  get an exact sequence of $T$-modules

$$0\to \begin{pmatrix}
K^1_1\\ K_2^1
\end{pmatrix}_{\varphi^{K^1}}\to\begin{pmatrix}
G^0_1\\(U\otimes_A G_1^0)\oplus G_2^0
\end{pmatrix}\stackrel{\begin{pmatrix}
f_1^0\\f_2^0
\end{pmatrix}}{\to} M\to 0.$$

Similarly, there exists an exact sequence of $B$-modules $G_2^1\stackrel{g_2^1}{\to} K_2^1\to 0$ where $G_2^1$ is $G_{C_2}$-projective and then, 
we  get an exact sequence of $T$-modules

$$0\to \begin{pmatrix}
K^2_1\\ K_2^2
\end{pmatrix}_{\varphi^{K^2}}\to\begin{pmatrix}
G^1_1\\(U\otimes_A G_1^1)\oplus G_2^1
\end{pmatrix}\to \begin{pmatrix}
K^1_1\\ K_2^1
\end{pmatrix}_{\varphi^{K^1}}\to 0.$$
repeat this process, we get the exact sequence of $T$-modules
$$0\to \begin{pmatrix}
0\\K^{m-k}_2
\end{pmatrix}\to \begin{pmatrix}
G^{m-k-1}_1\\(U\otimes_A G_1^{m-k-1})\oplus G_2^{m-k-1}
\end{pmatrix}
\stackrel{\begin{pmatrix}
f_1^{m-k-1}\\f_2^{m-k-1}
\end{pmatrix}}{\longrightarrow}$$

$$ \cdots\to\begin{pmatrix}
G^1_1\\(U\otimes_A G_1^1)\oplus G_2^1
\end{pmatrix}\stackrel{\begin{pmatrix}
f_1^1\\f_2^1
\end{pmatrix}}{\longrightarrow} \begin{pmatrix}
G^0_1\\(U\otimes_A G_1^0)\oplus G_2^0
\end{pmatrix}\stackrel{\begin{pmatrix}
f_1^0\\f_2^0
\end{pmatrix}}{\longrightarrow} M\to 0$$

Note that $\2Gpc((U\otimes_AG_1^i)\oplus G_2^i)=\2Gpc( U\otimes_AG_1^i)\leq k$, for every  $i\in\{0,\cdots,m-k-1\}$. So, by \cite[Proposition 3.11(2)]{BGO1} and the exact sequence 

$$0\to K^{m-k}_2
\to (U\otimes_A G_1^{m-k-1})\oplus G_2^{m-k-1}
\stackrel{f_2^{m-k-1}
}{\longrightarrow}\cdots\to (U\otimes_A G_1^0)\oplus G_2^0
\stackrel{f_2^0
}{\to} M_2\to 0$$

 we get that $\2Gpc(K_2^{m-k})\leq k$. This means that,  there exists an exact sequence of $B$-modules 
$$0\to G_2^m\to \cdots\to G_2^{m-k+1}\to G_2^{m-k}\to K_2^{m-k}\to 0.$$
Thus, There exists an exact sequence of $T$-modules
$$0\to \begin{pmatrix}
0\\G_2^m
\end{pmatrix}\to \cdots \to \begin{pmatrix}
0\\G_2^{m-k+1}
\end{pmatrix}\to$$
$$  \begin{pmatrix}
0\\G_2^{m-k}
\end{pmatrix}\to \begin{pmatrix}
G^{m-k-1}_1\\(U\otimes_A G_1^{m-k-1})\oplus G_2^{m-k-1}
\end{pmatrix}
\stackrel{\begin{pmatrix}
f_1^{m-k-1}\\f_2^{m-k-1}
\end{pmatrix}}{\longrightarrow}$$

$$\cdots \to \begin{pmatrix}
G^1_1\\(U\otimes_A G_1^1)\oplus G_2^1
\end{pmatrix}\stackrel{\begin{pmatrix}
f_1^1\\f_2^1
\end{pmatrix}}{\longrightarrow}\begin{pmatrix}
G^0_1\\(U\otimes_A G_1^0)\oplus G_2^0
\end{pmatrix}\stackrel{\begin{pmatrix}
f_1^0\\f_2^0
\end{pmatrix}}{\longrightarrow} M\to 0$$
By Theorem \ref{structure of $G_C$-projective}, all $\begin{pmatrix}
G^i_1\\(U\otimes_A G_1^i)\oplus G_2^i
\end{pmatrix}$ and all $\begin{pmatrix}
0\\G_2^j
\end{pmatrix}$ are $G_C$-projectives. Thus, $\Gpc(M)\leq m$. \cqfd

The following consequence of Theorem \ref{Estimation of dimension of a T-module} extends \cite[Proposition 2.8(1)]{EIT} and \cite[Theorem 2.7(1)]{ZLW} to the relative setting.
\begin{cor}Let $C=\textbf{p}(C_1,C_2)$ be w-tilting,  $U$ $C$-compatible and  $M=\begin{pmatrix}
M_1\\M_2
\end{pmatrix}_{\varphi^M}$ a $T$-module. If $SG_{C_2}-PD(B)<\infty,$
 then $\Gpc(M)<\infty$ if and only if $\1Gpc(M_1)<\infty$ and $\2Gpc(M_2)<\infty.$
\end{cor}

The following theorem gives an estimate of the left $G_C$-projective global dimension of $T$.

\begin{thm}\label{Estimation} Let $C=\textbf{p}(C_1,C_2)$ be w-tilting and $U$ $C$-compatible.
 Then 
$$max\{G_{C_1}-PD(A),G_{C_2}-PD(B)\}$$
$$\leq G_C-PD(T)\leq $$
$$ max\{G_{C_1}-PD(A)+SG_{C_2}-PD(B)+1,G_{C_2}-PD(B)\}.$$
\end{thm}
\proof We prove first that $max\{G_{C_1}-PD(A),G_{C_2}-PD(B)\}\leq G_C-PD(T)$. We may assume that $n:=G_C-PD(T)<\infty.$ Let $M_1$ be an $A$-module and $M_2$ be a $B$-module. Since $\Gpc(\begin{pmatrix}
M_1\\U\otimes_A M_2
\end{pmatrix}\leq n$ and  $\Gpc(\begin{pmatrix}
0\\M_2
\end{pmatrix}\leq n$, $\1Gpc(M_1)\leq n$ and $\2Gpc(M_2)\leq n$ by Lemma \ref{dimension of special modules}. Thus $G_{C_1}-PD(A)\leq n$ and 
$G_{C_2}-PD(B)\leq n$.

Next we prove that $$G_C-PD(T)\leq max\{G_{C_1}-PD(A)+1+ SG_{C_2}-PD(B),G_{C_2}-PD(B)\}.$$ We may assume that $$m:=max\{G_{C_1}-PD(A)+1+ SG_{C_2}-PD(B),G_{C_2}-PD(B)\}< \infty.$$ Then 
$n_1:=G_{C_1}-PD(A)<\infty$ and $k:=SG_{C_2}-PD(B)\leq n_2:=G_{C_2}-PD(B)<\infty$

Let $M=\begin{pmatrix}
M_1\\M_2
\end{pmatrix}_{\varphi^M}$ be a $T$-module. By Theorem  \ref{Estimation of dimension of a T-module}, $$\Gpc(M)\leq
  max\{n_1+k+1,n_2\}\leq m.$$
\cqfd

\begin{cor}\label{when rel gldim of T is finite} Let $C=\textbf{p}(C_1,C_2)$ be w-tilting and $U$ $C$-compatible. Then \\
	$ G_C-PD(T)<\infty $ if and only if $G_{C_1}-PD(A)<\infty$   and $G_{C_2}-PD(B)<\infty $
 
\end{cor}

Recall that a ring $R$ is called left Gorenstein regular if the category $R$-Mod is Gorenstein (\cite[Definition 2.1]{EIT} and \cite[Definition 2.18]{EEG}).

We know by \cite[Theorem 1.1]{BM}, that the following equality holds:
$$sup\{\Gpd_R(M)\;|\; \text{ $M\in R$-Mod}\}=sup\{\Gid_R(M)\;|\; \text{$M\in R$-Mod}\}.$$
and this common value is call the left global Gorenstein dimension of $R$, denoted by $l.\Ggldim(R)$. As a consequence of \cite[Theorem 2.28]{EEG}, a ring $R$ is left Gorenstein regular if and only if the global Gorenstein dimension of $R$ is finite.

We shall say that a ring $R$ is left $n$-Gorenstein regular if $n=l.\Ggldim(R)<\infty.$

Enochs, Izurdiaga and Torrecillas, characterized in  \cite[Theorem 3.1]{EIT} when $T$ is left Gorenstein regular under the conditions that $_BU$ has finite projective dimension and $U_A$ has finite flat dimension. As a direct consequence of Corollary \ref{when rel gldim of T is finite}, we refine this result.
\begin{cor}	Assume that $U$ is compatible. Then 
$T$ is left Gorenstein regular if and only if so are $A$ and $B$.

\end{cor}

There are some cases when the estimate in  Theorem \ref{Estimation} becomes an exact formula which computes left $G_C$-projective global dimension of $T$. 

Recall that an injective cogenerator $E$ in $R$-Mod is said to be strong if any $R$-module embeds in a direct sum of copies of $E$.

\begin{cor}  Let $C=\textbf{p}(C_1,C_2)$ be w-tilting and $U$ $C$-compatible.

\begin{enumerate}

\item If $U=0$ then 
$$ G_C-PD(T)= max\{G_{C_1}-PD(A),G_{C_2}-PD(B)\}$$

\item If $A$ is left noetherian and $_AC_1$ is a strong injective cogenerator, then 

$$G_C-PD(T)=\begin{cases}
G_{C_2}-PD(B) & \text{if $U=0$ }\\
max\{SG_{C_2}-PD(B)+1,G_{C_2}-PD(B)\} & \text{if  $U\neq 0$}
\end{cases}$$
\end{enumerate}

\end{cor} 
\proof 1.  Using a similar way as we do in the proof of Theorem \ref{Estimation of dimension of a T-module} and \ref{Estimation}, we can prove this statment. We only need to notice that if $U=0$, then  a $T$-module $M=\begin{pmatrix}
M_1\\M_2
\end{pmatrix}_{\varphi^M}$ is $G_C$-projective if and only if $M_1$ is $G_{C_1}$-projective and $M_2$ is $G_{C_2}$-projective (since $\varphi^M$ is always injective and $M_2=\overline{M}_2$) by Theorem \ref{structure of $G_C$-projective}.

2. Note first that $G_{C_1}-PD(A)=0$ by  \cite[Corollary 2.3]{BGO2}. Then the case $U=0$ follows by 1. Assume that $U\neq 0$. Note that  by Theorem \ref{structure of $G_C$-projective}, $ \begin{pmatrix}
A\\0
\end{pmatrix}$ is not $G_C$-projective since $U\neq 0$. Hence  $G_{C_2}-PD(B)\geq \Gpc_T( \begin{pmatrix}
A\\0
\end{pmatrix})\geq 1$.

By Theorem \ref{Estimation}, we have the inequality  $$G_{C_2}-PD(B)\leq G_C-PD(T)\leq 
 max\{SG_{C_2}-PD(B)+1,G_{C_2}-PD(B)\}.$$
So, the case $SG_{C_2}-PD(B)+1\leq G_{C_2}-PD(B)$ is clear and  we only need to prove the result when  $SG_{C_2}-PD(B)+1> n:=G_{C_2}-PD(B)$. Since $\2Gpc(U\otimes_A G)\leq G_{C_2}-PD(B)=n$ for every  $G\in G_{C_1}P(A)$, $SG_{C_2}-PD(B)=n$. Let $G_1$ be a  $G_{C_1}$-projective $A$-module with $\2Gpc(U\otimes_A G_1)=n$ and consider the following short exact sequnece 
$$0\to \begin{pmatrix}
0\\U\otimes_A G_1
\end{pmatrix}\to \begin{pmatrix}
G_1\\U\otimes_A G_1
\end{pmatrix}\to \begin{pmatrix}
G_1\\0
\end{pmatrix}\to 0.$$
By Theorem  \ref{structure of $G_C$-projective}, $ \begin{pmatrix}
G_1\\U\otimes_A G_1
\end{pmatrix}$ is $G_C$-projective and  by Lemma \ref{dimension of special modules} $$\Gpc( \begin{pmatrix}
0\\U\otimes_A G_1
\end{pmatrix})=\2Gpc(U\otimes_A G)=n.$$ Thus by \cite[Proposition 3.11(4)
]{BGO1}
$$\Gpc( \begin{pmatrix}
G_1\\0 
\end{pmatrix})=\Gpc( \begin{pmatrix}
0\\U\otimes_A G_1
\end{pmatrix})+1=n+1=SG_{C_2}-PD(B)+1.$$
This shows that $G_C-PD(T)=SG_{C_2}-PD(B)+1$ and the proof is finished. \cqfd

\begin{cor}  Let $R$ be a ring, $T(R)=\begin{pmatrix}
R&0\\R&R
\end{pmatrix}$ and $C=\textbf{p}(C_1,C_1)$ where $C_1$ is w-tilting. Then

$$ G_{C}-PD(T(R))
=G_{C_1}-PD(R)+1.$$

\end{cor}
\proof
 Note first that $C$ is w-tilting $T(R)$-module, $R$ is $C$-compatible and $SG_{C_1}-PD(R)=0$. So, by Theorem \ref{Estimation}, 
$$G_{C_1}-PD(R)\leq G_{C}-PD(T(R))
\leq G_{C_1}-PD(R)+1.$$
The case $G_{C_1}-PD(R)=\infty$ is clear. Assume that $n:=G_{C_1}-PD(R)<\infty.$ 

There exists an $R$-module $M$ with $\1Gpc(M)=n$ and    $\Ext_R^n(M,X)\neq 0$ for some $X\in \Add_R(C_1)$ by \cite[Theorem 3.8]{BGO1}. If we apply the functor $\Hom_{T(R)}(-,\begin{pmatrix}
0\\X
\end{pmatrix})$ to the exact sequence of $T(R)$-modules
$$0\to \begin{pmatrix}
0\\M
\end{pmatrix}\to \begin{pmatrix}
M\\M
\end{pmatrix}_{1_M}\to \begin{pmatrix}
M\\0
\end{pmatrix}\to 0$$
we get an exact sequence 
$$\cdots \to\Ext^n_{T(R)}(\begin{pmatrix}
M\\M
\end{pmatrix},\begin{pmatrix}
0\\X
\end{pmatrix})\to \Ext^n_{T(R)}(\begin{pmatrix}
0\\M
\end{pmatrix},\begin{pmatrix}
0\\X
\end{pmatrix})\to $$$$ \Ext^{n+1}_{T(R)}(\begin{pmatrix}
M\\0
\end{pmatrix},\begin{pmatrix}
0\\X
\end{pmatrix}) \to \Ext^{n+1}_{T(R)}(\begin{pmatrix}
M\\M
\end{pmatrix},\begin{pmatrix}
0\\X
\end{pmatrix})\to\cdots $$
By Lemma \ref{Ext}, $\Ext^{i\geq 1}_{T(R)}(\begin{pmatrix}
M\\M
\end{pmatrix},\begin{pmatrix}
0\\X
\end{pmatrix})\cong \Ext_R^{i\geq 1}(M,0)=0$. Again by Lemma \ref{Ext} and the above exact sequence, $$\Ext^{n+1}_{T(R)}(\begin{pmatrix}
M\\0
\end{pmatrix},\begin{pmatrix}
0\\X
\end{pmatrix})\cong \Ext^n_{T(R)}(\begin{pmatrix}
0\\M
\end{pmatrix},\begin{pmatrix}
0\\X
\end{pmatrix})\cong \Ext_R^n(M,X)\neq 0.$$ 
 since $\begin{pmatrix}
0\\X
\end{pmatrix}\in \Add_{T(R)}(C)$ by Lemma \ref{Add-Prod}(1), it follows that $n<\Gpc(\begin{pmatrix}
M\\0
\end{pmatrix})$ by \cite[Theorem 3.8]{BGO1}. But $\Gpc(\begin{pmatrix}
M\\0
\end{pmatrix})\leq G_C-PD(T(R))
\leq n+1$. Thus $\Gpc(\begin{pmatrix}
M\\0
\end{pmatrix})=n+1$. which means that $G_C-PD(T(R))
= n+1.$ \cqfd
\begin{cor} Let $R$ be a ring, $T(R)=\begin{pmatrix}
	R&0\\R&R
	\end{pmatrix}$ and $n\geq 0$ an integer. Then  $T(R)$ is left $(n+1)$-Gorenstein regular if and only if $R$ left $n$-Gorenstein regular .

\end{cor}

The authors in \cite{AS} establish a relationship between the projective  dimension of
modules over $T$ and modules over $A$ and $B$. Given  an integer $n\geq 0$  and $M=\begin{pmatrix}
M_1\\M_2
\end{pmatrix}_{\varphi^M}$  a $T$-module, they proved that $\pd_T(M)\leq n $ if and only if $\pd_A(M_1)\leq n$, $\pd_B(\overline{M}_2)\leq n$ and the map related to the $n$-th syzygy of $M$ is injective. The following example shows that this is not true in general.

\begin{exmp}\label{counter example} Let $R$ be a left hereditary ring which not semisimple and  let $T(R)=\begin{pmatrix}
R&0\\R&R
\end{pmatrix}$. Then $lD(T(R))=lD(R)+1=2 $ by \cite[corollary 3.4(3)]{M}. This means that there exists a $T(R)$-module $M=\begin{pmatrix}
M_1\\M_2
\end{pmatrix}_{\varphi^M}$ with $\pd_{T(R)}(M)=2$. If $K^1=\begin{pmatrix}
K_1^1\\K_2^1
\end{pmatrix}_{\varphi^{K^1}}$ is the first  syzygy of $M$, then there exists an exact sequence of $T(R)$-modules $$0\to K^1\to P\to M\to 0$$
where $P=\begin{pmatrix}
P_1\\P_2
\end{pmatrix}_{\varphi^P}$ is projective. Then we get the following commutative diagram

$$   \xymatrix{
 &  &0\ar[d]  &   & &     \\
 0\ar[r]& K^1_1\ar[d]_{\varphi^{K^11}} \ar[r]  &  P_1\ar[d]_{\varphi^P} \ar[r]& M_1\ar[d]_{\varphi^M} \ar[r] &  0   \\
 0\ar[r]& K^1_2\ar[d] \ar[r] & P_2\ar[d] \ar[r] &  M_2 \ar[d]  \ar[r] & 0\\
& \overline{K^1_2}\ar[r] \ar[d] &  \overline{P}_2  \ar[d] \ar[r]& \overline{M}_2\ar[r] \ar[d] & 0\\
&0  &0 &0  &    }  $$

By Snake Lemma $\varphi^{K^1}$ is injective. On the other hand, Since $lD(R)=1$, $\pd_R(M_1)\leq 1$ and $\pd_R(\overline{M}_2)\leq 1$. But $\pd_{T(R)}(M)=2>1$. \cqfd

\end{exmp}

\noindent {\bf Acknowledgement.}
The third and forth authors were partially supported by Ministerio de
Econom  \'{\i}a y Competitividad, grant reference 2017MTM2017-86987-P, and Junta de
Andaluc\'{\i}a, grant reference
P20-00770. The authors would like to thank Professor Javad Asadollahi for the discussion on the Example  \ref{counter example}.

\end{document}